\documentclass[10pt,reqno]{amsart}
\usepackage[english]{babel}
\usepackage[utf8]{inputenc}
\usepackage{amssymb,amsmath,amsfonts,amsthm,amscd,latexsym,graphicx,epsfig,colordvi}
\usepackage{comment}

\usepackage{geometry} 
\let\oldsection\section
\renewcommand{\section}{
  \renewcommand{\theequation}{\thesection.\arabic{equation}}
  \oldsection}
\let\oldsubsection\subsection
\renewcommand{\subsection}{
  \renewcommand{\theequation}{\thesubsection.\arabic{equation}}
  \oldsubsection}

\begin{document}

\author{A.S.Panasenko}
\address{Alexander Sergeevich Panasenko
\newline\hphantom{iii} Novosibirsk State University,
\newline\hphantom{iii} 1, Pirogova str.
\newline\hphantom{iii} 630090, Novosibirsk, Russia
\newline\hphantom{iii} Sobolev Institute of Mathematics,
\newline\hphantom{iii} pr. Koptyuga, 4,
\newline\hphantom{iii} 630090, Novosibirsk, Russia
}
\email{a.panasenko@g.nsu.ru}%


\title{ROTA-BAXTER OPERATORS OF NONZERO WEIGHT ON THE SPLIT OCTONIONS}
\maketitle
\small
{\begin{quote}
\noindent{\sc Abstract. }  We describe Rota-Baxter operators on split octonions. It turns out that up to some transformations there exists exactly one such non-splitting operator over any field. We also obtain a description of all decompositions of split octonions over a quadratically closed field of characteristic different from 2 into a sum of two subalgebras, which describes the splitting Rota-Baxter operators. It completes the classification of Rota-Baxter operators on composition algebras of any weight.

\medskip

\end{quote}}

\section{Introduction}

The notion of a composition algebra arose from the natural question of representing the product of a sum of $n$ squares and a sum of $n$ squares again as a sum of two squares, the so-called Hurwitz problem. Although every composition algebra has dimension $1$, $2$, $4$, or $8$, composition division algebras may differ depending on the underlying field. For example, over the field of real numbers there are exactly 4 composition division algebras, but over the algebraic number field there may be arbitrarily many nonisomorphic 4-dimensional and 8-dimensional composition division algebras (\cite{Springer}).

The situation is quite different for composition algebras with zero divisors, the so-called split algebras. It is well known that over any field of characteristic not 2 there exist exactly 4 such algebras: one in each dimension 1, 2, 4, 8. The two-dimensional algebra is isomorphic to the direct sum of two copies of the ground field, and the four-dimensional algebra is isomorphic to the algebra of second-order matrices over the ground field. Thus, the study of 8-dimensional composition algebras with zero divisors (split octonions, split Cayley-Dickson algebra) is a natural generalization of the study of the direct sum of fields and second-order matrices.

Rota-Baxter operators arise naturally from several areas. The Rota-Baxter operator is a generalization of the integration operator \cite{Baxter}. In the paper \cite{Semenov}, Rota-Baxter operators of nonzero weight appeared as solutions of the modified Yang-Baxter equation. There are a lot of applications of Rota-Baxter operators to various areas of algebra. In \cite{Aguiar} and \cite{An}, the connection of Rota-Baxter operators with pre-Lie algebras is investigated. In \cite{Double}, the connection of Rota-Baxter operators with double Lie algebras is investigated. Moreover, to any decomposition of an algebra into a direct sum of two subalgebras there corresponds some Rota-Baxter operator of nonzero weight. Thus, the study of Rota-Baxter operators on useful algebras is of interest because of the importance of studying the structure of subalgebras. 

Descriptions of Rota-Baxter operators of nonzero weight are known on a lot of simple algebras and superalgebras: a~simple Jordan algebra of a~bilinear form of odd dimension, $M_2(F)$, $K_3$ (all in \cite{BGP}), $M_3(\mathbb{C})$ \cite{M1,M2,M3}.

One of the most important works in the direction of studying the structure of subalgebras in the octonions is the paper by Knarr and Stroppel \cite{Octonions}. In it, the authors described all subalgebras up to an automorphism. Moreover, it was proved that two subalgebras lie in the same orbit with respect to the automorphism group of the octonions if and only if they are isomorphic as algebras, with the exception of the case of one-dimensional idempotent algebras.

It is well known that over a fixed field $F$ of characteristic different from 2, there exist only three split composition algebras: the direct sum of two fields $F$, the matrix algebra $M_2(F)$, and the split octonions over $F$. It follows from the results of \cite{BGP} that all Rota-Baxter operators on non-split composition algebras are trivial. In \cite{An}, Rota-Baxter operators on the direct sum of two fields were described (later the description was generalized to the direct sum of a finite number of fields, \cite{GubarevField}). In \cite{BGP}, Rota-Baxter operators on the algebra of second-order matrices over an algebraically closed field were described. Thus, the problem of classifying Rota-Baxter operators on composition algebras was reduced to the problem of describing Rota-Baxter operators on split octonions. In the paper \cite{RBZ} Rota-Baxter operators of zero weight on split octonions were described.

In this paper we describe Rota-Baxter operators on split octonions. It turns out that up to some transformations there exists exactly one such non-splitting operator over any field. We also obtain a description of all decompositions of split octonions over a quadratically closed field of characteristic $\neq 2$ into a sum of two subalgebras, which classifies the splitting Rota-Baxter operators. Thus, this paper completes the classification of Rota-Baxter operators on composition algebras.

\section{Preliminaries}

\textbf{Definition.} Let $A$ be an algebra over a field $F$. A linear map $R: A \to A$ is called a \textit{\textbf{Rota-Baxter operator}} of weight $\lambda\in F$ if for any $x,y\in A$ we have
\[R(x)R(y) = R(R(x)y+xR(y)+\lambda xy).\]

Let $A = A_1 \oplus A_2$ be a decomposition of $A$ into a direct sum of two subalgebras $A_1$ and $A_2$. Then a linear transformation $R$ on $A$, defined on the subalgebras as follows: $R|_{A_1}=-\lambda id$, $R|_{A_2}=0$, is a Rota-Baxter operator of weight $\lambda$ on $A$.

\smallskip\textbf{Definition.} Let $R$ be a Rota-Baxter operator of nonzero weight $\lambda$ on an algebra $A$. The operator $R$ is called \textbf{\textit{splitting}} if there exists a decomposition $A= A_1\oplus A_2$ into a direct sum of two subalgebras, where $R|_{A_1} = -\lambda id$, $R|_{A_2}=0$.

In particular, there is a natural bijection between the set of splitting Rota-Baxter operators of unit weight and decompositions of $A$ into a direct sum of two subalgebras.

Note that a Rota-Baxter operator $R$ of nonzero weight $\lambda$ is splitting if and only if $R^2 = -\lambda R$.

In this paper we consider a split Cayley-Dickson algebra $\mathbb{O} = M_2(F) + vM_2(F)$ with the following multiplication table:
\[a\cdot b = ab, \quad a\cdot vb = v(\overline{a}b), \quad va \cdot b = v(ba), \quad va\cdot vb = b\overline{a},\]
where $x\cdot y$ is a multiplication in $\mathbb{O}$, $ab$ is a multiplication in $M_2(F)$, $\overline{a}$ is a symplectic involution on $M_2(F)$:
\[\overline{\begin{pmatrix}
a & b \\
c & d
\end{pmatrix}} = \begin{pmatrix}
d & -b \\
-c & a
\end{pmatrix}.\]

\smallskip In what follows we will need the following examples of automorphisms and antiautomorphisms on a split Cayley-Dickson algebra (see~\cite{RBZ}).

\medskip\textbf{Example 1.} \textit{Let $\varphi:\mathbb{O}\to \mathbb{O}$ be a linear map such that $\varphi(e_{ij})=e_{ij}$ for any $i,j\in\{1,2\}$, $\varphi(ve_{12})=ve_{12}$, $\varphi(ve_{22})=ve_{22}$, $\varphi(ve_{11})=ve_{11}+\alpha ve_{12}$, $\varphi(ve_{21})=ve_{21}+\alpha ve_{22}$ for some $\alpha\in F$. Then $\varphi$ is an automorphism on the algebra $\mathbb{O}$.}

\medskip\textbf{Example 2.} \textit{Let $\varphi:\mathbb{O}\to \mathbb{O}$ be a linear map such that $\varphi(e_{11}) = e_{11} + \alpha ve_{22}$, $\varphi(e_{12}) = e_{12} + \alpha ve_{12}$, $\varphi(e_{21}) = e_{21}$, $\varphi(e_{22}) = e_{22} - \alpha ve_{22}$, $\varphi(ve_{11}) = ve_{11}-\alpha e_{11}+\alpha e_{22}-\alpha^2 ve_{22}$, $\varphi(ve_{12}) = ve_{12}$, $\varphi(ve_{21}) = ve_{21} + \alpha e_{21}$, $\varphi(ve_{22}) = ve_{22}$ for some $\alpha\in F$. Then $\varphi$ is an automorphism on the algebra $\mathbb{O}$. }

\medskip\textbf{Example 3.} \textit{Let $\varphi:\mathbb{O}\to \mathbb{O}$ be a linear map such that $\varphi(e_{ii}) = e_{ii}$, $\varphi(ve_{ii}) = ve_{ii}$ for $i\in\{1,2\}$, $\varphi(e_{12}) = e_{12}$, $\varphi(ve_{12}) = ve_{12}$, $\varphi(e_{21}) = e_{21}+\alpha ve_{12}$, $\varphi(ve_{21}) = ve_{21} + \alpha e_{12}$ for some $\alpha\in F$. Then $\varphi$ is an automorphism on the algebra~$\mathbb{O}$.}

\medskip\textbf{Example 4.} \textit{Let $\varphi:\mathbb{O}\to \mathbb{O}$ be a linear map such that $\varphi(e_{ii}) = e_{ii}$, $\varphi(e_{12}) = e_{12}$, $\varphi(e_{21}) = e_{21} + \alpha ve_{11}$, $\varphi(e_{22}) = e_{22}$, $\varphi(ve_{11}) = ve_{11}$, $\varphi(ve_{12}) = ve_{12}$, $\varphi(ve_{21}) = ve_{21}$, $\varphi(ve_{22}) = ve_{22} - \alpha e_{12}$ for some $\alpha\in F$. Then $\varphi$ is an automorphism on the algebra $\mathbb{O}$.}

\medskip\textbf{Example 5.} \textit{Let $\varphi:\mathbb{O}\to \mathbb{O}$ be a linear map such that $\varphi(e_{11}) = e_{11} - \alpha ve_{12}$, $\varphi(e_{12}) = e_{12}$, $\varphi(e_{21}) = e_{21} + \alpha ve_{22}$, $\varphi(e_{22}) = e_{22} + \alpha ve_{12}$, $R(ve_{11}) = ve_{11} - \alpha e_{12}$, $R(ve_{12}) = ve_{12}$, $\varphi(ve_{21}) = ve_{21}-\alpha e_{11}+\alpha e_{22}+\alpha^2 ve_{12}$, $\varphi(ve_{22}) = ve_{22}$ for some $\alpha\in F$. Then $\varphi$ is an automorphism on the algebra $\mathbb{O}$. }

\medskip\textbf{Example 6.} \textit{Let $\varphi:\mathbb{O}\to \mathbb{O}$ be a linear map such that $\varphi(e_{11}) = e_{11} + \alpha e_{12}$, $\varphi(e_{12}) = e_{12}$, $\varphi(e_{21}) = e_{21} - \alpha e_{11} + \alpha e_{22} - \alpha^2 e_{12}$, $\varphi(e_{22}) = e_{22} - \alpha e_{12}$, $R(ve_{11}) = ve_{11}$, $R(ve_{12}) = ve_{12}$, $\varphi(ve_{21}) = ve_{21}-\alpha ve_{11}$, $\varphi(ve_{22}) = ve_{22} - \alpha ve_{12}$ for some $\alpha\in F$. Then $\varphi$ is an automorphism on the algebra $\mathbb{O}$. }

\medskip\textbf{Example 7.} \textit{Let $\varphi:\mathbb{O}\to \mathbb{O}$ be a linear map such that $\varphi(e_{ij})=e_{ij}$, for any $i,j\in\{1,2\}$, $\varphi(ve_{11}) = \alpha ve_{11}$, $\varphi(ve_{21})=\alpha ve_{21}$, $\varphi(ve_{22}) = \frac{1}{\alpha} ve_{22}$, $\varphi(ve_{12})= \frac{1}{\alpha} ve_{12}$ for some $0\neq \alpha\in F$. Then $\varphi$ is an automorphism on the algebra $\mathbb{O}$.}

\medskip\textbf{Example 8.} \textit{Let $\varphi:\mathbb{O}\to \mathbb{O}$ be a~linear map such that $\varphi(e_{ii}) = e_{ii }$, $\varphi(e_{12}) = e_{12} - \alpha ve_{22}$, $\varphi(e_{21}) = e_{21}$, $\varphi(e_{22}) = e_{22}$, $\varphi(ve_{11}) = ve_{11} + \alpha e_{21}$, $\varphi( ve_{12}) = ve_{12}$, $\varphi(ve_{21}) = ve_{21}$, $\varphi(ve_{22}) = ve_{22}$ for some $\alpha\in F$. Then $\varphi$ is an~automorphism on the algebra $\mathbb{O}$.}

\medskip\textbf{Example 9.} \textit{Let $\varphi:\mathbb{O}\to \mathbb{O}$ be a~linear map such that $\varphi(e_{ij})=e_{ij }$ for any $i,j\in\{1,2\}$, $\varphi(ve_{12})=-ve_{11}$, $\varphi(ve_{11})=ve_{12}$, $\varphi(ve_{21})=ve_{22}$, $\varphi(ve_{22})=-ve_{21}$. Then $\varphi$ is an~automorphism on the algebra $\mathbb{O}$.}

\medskip\textbf{Example 10.} \textit{Let $\varphi:\mathbb{O}\to \mathbb{O}$ be a linear map such that $\varphi(e_{ii}) = e_{ii} $, $\varphi(ve_{ii}) = ve_{ii}$ for $i\in\{1,2\}$, $\varphi(e_{21}) = e_{21}$, $\varphi(ve_{21}) = ve_{21}$, $\varphi(e_{12}) = e_{12}+\alpha ve_ {21}$, $\varphi(ve_{12}) = ve_{12} + \alpha e_{21}$ for some $\alpha\in F$. Then $\varphi$ is an automorphism on the algebra~$\mathbb{O}$.}

\medskip\textbf{Example 11.} \textit{Let $\varphi:\mathbb{O}\to \mathbb{O}$ be a~linear map such that $\varphi(e_{ii}) = e_{ii }$ for $i\in\{1,2\}$, $\varphi(e_{12}) = e_{21}$, $\varphi(e_{21}) = e_{12}$, $\varphi(ve_{11}) = - ve_{21}$, $\varphi(ve_{12}) = - ve_{22}$, $\varphi(ve_{21}) = ve_{11}$, $\varphi(ve_{22}) = ve_{12}$ for some $\alpha\in F$. Then $\varphi$ is an~antiautomorphism on the algebra $\mathbb{O}$.}

\medskip\textbf{Example 12.} \textit{Let $\varphi:\mathbb{O}\to \mathbb{O}$ be a~linear map such that $\varphi(e_{ii})=e_{ii }$, for any $i\in\{1,2\}$, $\varphi(e_{12}) = \alpha e_{12}$, $\varphi(ve_{11}) = \alpha ve_{11}$, $\varphi(e_{21})=\frac{1}{\alpha} e_{21}$, $\varphi(ve_{22 }) = \frac{1}{\alpha} ve_{22}$, $\varphi(ve_{12})=ve_{12}$, $\varphi(ve_{21})=ve_{21}$ for some $0\neq \alpha\in F$. Then $\varphi$ is an~automorphism on the algebra $\mathbb{O}$.}

\medskip\textbf{Example 13.} \textit{Let $\varphi:\mathbb{O}\to \mathbb{O}$ be a~linear map such that $\varphi(e_{ii})=e_{ii }$, $\varphi(ve_{ii})=ve_{ii}$ for any $i\in\{1,2\}$, $\varphi(e_{12}) = \alpha e_{12}$, $\varphi(ve_{12}) = \alpha ve_{12}$, $\varphi(e_{21})=\frac{1 }{\alpha} e_{21}$, $\varphi(ve_{21}) = \frac{1}{\alpha} ve_{21}$ for some $0\neq \alpha\in F$. Then $\varphi$ is an~automorphism on the algebra $\mathbb{O}$.}

\medskip\textbf{Example 14.} \textit{Let $\varphi:\mathbb{O}\to \mathbb{O}$ be a~linear map such that $\varphi(e_{11})=e_{22 }$, $\varphi(e_{12})=ve_{22}$, $\varphi(e_{21})=ve_{11}$, $\varphi(e_{22})=e_{11}$ and $\varphi^2 =\mathrm{id}$. Then $\varphi$ is an~involution on the algebra $\mathbb{O}$.}

\medskip\textbf{Example 15.} \textit{Let $\varphi:\mathbb{O}\to \mathbb{O}$ be a~linear map such that $\varphi(e_{ij})=e_{ij}$ for any $i,j\in\{1,2\}$, $\varphi(ve_{11})=ve_{11}$, $\varphi(ve_{21})=ve_{21}$, $\varphi(ve_{12})=ve_{12}+\alpha ve_{11}$, $\varphi(ve_{22})=ve_{22}+\alpha ve_{21}$ for some $\alpha\in F$. Then $\varphi$ is an~automorphism on the algebra $\mathbb{O}$.}

\medskip\textbf{Example 16.} \textit{Let $\varphi:\mathbb{O}\to \mathbb{O}$ be a~linear map such that $\varphi(e_{11}) = e_{11} + \alpha ve_{22}$, $\varphi(e_{12}) = e_{12} + \alpha ve_{12}$, $\varphi(e_{21}) = e_{21}$, $\varphi(e_{22}) = e_{22} - \alpha ve_{22}$, $\varphi(ve_{11}) = ve_{11}-\alpha e_{11}+\alpha e_{22}-\alpha^2 ve_{22}$, $\varphi(ve_{12}) = ve_{12}$, $\varphi(ve_{21}) = ve_{21} + \alpha e_{21}$, $\varphi(ve_{22}) = ve_{22}$ for some $\alpha\in F$. Then $\varphi$ is an~automorphism on the algebra $\mathbb{O}$.
}

\medskip\textbf{Example 17.} \textit{Let $\varphi:\mathbb{O}\to \mathbb{O}$ be a~linear map such that $\varphi(e_{11}) = e_{21}+e_{22}+ve_{12}$, $\varphi(e_{12}) = -e_{11}-e_{12}+e_{21}+e_{22}+ve_{12}$, $\varphi(e_{21}) = -e_{21}-ve_{22}$, $\varphi(e_{22}) = e_{11}-e_{21}-ve_{12}$, $\varphi(ve_{11}) = -e_{12}+ve_{11}-ve_{21}$, $\varphi(ve_{12}) = ve_{12}-ve_{22}$, $\varphi(ve_{21}) = -e_{11}+e_{21}+e_{22}+ve_{12}+ve_{21}$, $\varphi(ve_{22}) = ve_{22}$. Then $\varphi$ is an~antiautomorphism on the algebra $\mathbb{O}$.}

\medskip\textbf{Example 18.} \textit{Let $\varphi:\mathbb{O}\to \mathbb{O}$ be a~linear map such that $\varphi(e_{11}) = e_{11} + ve_{12}$, $\varphi(e_{12}) = e_{12}$, $\varphi(e_{21}) = e_{21}-ve_{22}$, $\varphi(e_{22}) = e_{22}-ve_{12}$, $\varphi(ve_{11}) = ve_{11}+e_{12}+ve_{12}$, $\varphi(ve_{12}) = -e_{12}-ve_{11}$, $\varphi(ve_{21}) = ve_{21}+e_{11}-e_{22}+ve_{12}+ve_{22}$, $\varphi(ve_{22}) = -e_{11}+e_{22}-ve_{12}-ve_{21}$. Then $\varphi$ is an~automorphism on the algebra $\mathbb{O}$.}

\section{Non-splitting Rota-Baxter operators of non-zero weight with non-three-dimensional kernel}

We will need the following well-known result.

\medskip\textbf{Proposition 1 (\cite{BGP}).} \textit{Let $R$ be a Rota-Baxter operator of nonzero weight on a finite-dimensional unital algebra over $F$, with $R(1)\in F$. Then $R$ is splitting.}

\medskip First, we prove that non-splitting Rota-Baxter operators on octonions can have only three-dimensional or four-dimensional kernels.

\medskip\textbf{Lemma 1.} \textit{Let $R$ be a Rota-Baxter operator of nonzero weight on a split Cayley-Dickson algebra $\mathbb{O}$ over a field $F$ and $\dim(\mathrm{Ker}(R)) \notin\{3,4\}$. Then $R$ is splitting.}

\textbf{Proof.} Recall (see, e.g., \cite{BGP}) that $\mathrm{Ker}(R)$ and $\mathrm{Im}(R)$ are subalgebras of $\mathbb{O}$. Denote $n=\dim\mathrm{Ker}(R)$. By the results of \cite{Octonions}, there is no subalgebra of dimension~7 in $\mathbb{O}$, so $n\neq 7$.

If $n\in\{5,6,8\}$, then by the results of \cite{Octonions}, the subalgebra $\mathrm{Ker}(R)$ contains the unit of $\mathbb{O}$. It means that $R(1)=0$. Proposition~1 implies that $R$ is splitting.

If $n\le 2$, then $\dim\mathrm{Im}(R)\ge 6$. We introduce the notation: \[\mathbb{O}_0 = \{x\in \mathbb{O}\mid t(x)=0\} = vM_2(F) + Fe_{12} + Fe_{21} + F(e_{11}-e_{22}).\]
Then $\dim(R(\mathbb{O}_0))\ge 5$. It is well known (see, e.g., \cite{Springer}, section 1.8) that if $V$ is a subspace of $\mathbb{O}$ such that $n(x) = 0$ for any $x\in V$, then $\dim V\le 4$. Thus, there exists $y\in R(\mathbb{O}_0)$ such that $n(y)\neq 0$. By (\cite{BGP}, Lemma 3.1.c) the operator $R$ is splitting. The lemma is proved.

\medskip The following example of a non-splitting Rota-Baxter operator is well known.

\medskip\textbf{Example 19 (\cite{Spectrum}, Proposition~3).} {\it Let $\mathbb{O}$ be a split Cayley-Dickson algebra and $R_1:\mathbb{O}\to\mathbb{O}$ be a linear map on $\mathbb{O}$ defined as follows:
\[R_1(e_{22}) = e_{11},\quad R_1(e_{21}) = -e_{21},\quad R_1(ve_{21}) = -ve_{21},\quad R_1(ve_{22}) = -ve_{22}, \quad R_1(e_{1i})=R_1(ve_{1i})=0\]
for $i\in\{1,2\}$. Then $R_1$ is a non-splitting Rota-Baxter operator on $\mathbb{O}$.}

\medskip It turns out that $R_1$ is the only non-splitting operator on $\mathbb{O}$ with a four-dimensional kernel.

\medskip\textbf{Lemma 2.} \textit{Let $R$ be a non-splitting Rota-Baxter operator of unit weight on a split Cayley-Dickson algebra $\mathbb{O}$ and $\mathrm{dim}(\mathrm{Im}(R)) = 4$. Then, up to conjugation by an automorphism and anti-automorphism, the operator $R$ coincides with $R_1$.
}

\textbf{Proof.} Since $R$ is non-splitting by Proposition~1 and by the results of \cite{Octonions} we can assume that, up to automorphism and antiautomorphism, $\mathrm{Ker}(R) = Fe_{11}+Fe_{12}+Fve_{11}+Fve_{12}$. Note that in this case, it follows from the theorem on the sum of the dimensions of the kernel and the image that the elements $R(e_{21}),R(e_{22}),R(ve_{21}),R(ve_{22})$ form a basis of $\mathrm{Im}(R)$.

Let $R(e_{22}) = \sum\alpha_{ij}e_{ij} + \sum\beta_{ij}ve_{ij}$. Then
\[0 = R(e_{22})R(e_{12}) = R(R(e_{22})e_{12}) = \alpha_{21}R(e_{22}),\]
whence $\alpha_{21} = 0$. Next,
\[0 = R(ve_{12})R(e_{22}) = R(ve_{12}R(e_{22})) = -\beta_{21}R(-e_{22}),\]
whence $\beta_{21} = 0$. Next,
\[0 = R(ve_{11})R(e_{22}) = R(ve_{11}R(e_{22})) = \beta_{22}R(e_{22}),\]
whence $\beta_{22} = 0$. Next,
\begin{equation*}
\alpha_{11}^2e_{11} + (\alpha_{11}\alpha_{12}+\alpha_{12}\alpha_{22})e_{12} + \alpha_{22}^2 e_{22} + (\alpha_{22}\beta_{11}+\alpha_{11}\beta_{11})ve_{11} + (\alpha_{22}\beta_{12}+\beta_{12}\alpha_{11} )ve_{12} = R(e_{22})R(e_{22}) = (2\alpha_{22}+1)R(e_{22}),
\end{equation*}
where $\alpha_{22}^2 = (2\alpha_{22}+1)\alpha_{22}$, i.e. $\alpha_{22}\in\{0,-1\}$.

Suppose that $\alpha_{11}-\alpha_{22} \neq 1$. Then, comparing the coefficients of $e_{12},ve_{11},ve_{12}$ in the formula for $R(e_{22})^2$, we obtain $\alpha_{12} = \beta_{11} = \beta_{12} = 0$, that is, $R(e_{22}) = \alpha_{11}e_{11} + \alpha_{22}e_{22}$, and $\alpha_{11}^2 = (2\alpha_{22}+1)\alpha_{11}$. If $\alpha_{11} = 0$, then $\alpha_{22} = -1$ (otherwise $\alpha_{22} = 0$ and $R(1) = R(e_{22}) = 0$, so the operator $R$ is splitting by Proposition~1). But then $\alpha_{11}-\alpha_{22}=1$, it is a contradiction. Thus, $\alpha_{11}\neq 0$ and $\alpha_{11} = 2\alpha_{22} + 1$. If $\alpha_{22} = 0$, then $\alpha_{11} = 1$, it is a contradiction (in this case $\alpha_{11}-\alpha_{22} = 1$). Thus, $\alpha_{22} = -1$ and $\alpha_{11} = -1$. It means that $R(1)=R(e_{22}) = -1$ and the operator $R$ is splitting by Proposition~1, it is a contradiction.

So, $\alpha_{11}-\alpha_{22} = 1$. If $\alpha_{11} = 0$, then $\alpha_{22} = -1$ and $R(e_{22}) = -e_{22} + \alpha_{12} + \beta_{11}ve_{11} + \beta_{12}ve_{12}$. If $\alpha_{11}\neq 0$, then $\alpha_{11} = 2\alpha_{22}+1$, whence $\alpha_{22} = 0$, $\alpha_{11} = 1$, so that $R(e_{22}) = e_{11} + \alpha_{12}e_{12} + \beta_{11}ve_{11} + \beta_{12}ve_{12}$. In any case, if $\beta_{11}\neq 0$, then the automorphism from Example~2 with scalar $\alpha = -\frac{\beta_{12}}{\beta_{11}}$ allows us to assume that $\beta_{12} = 0$. Thus, we can assume that either $\beta_{11} = 0$ or $\beta_{12} = 0$. Then the automorphism from Example~4 allows us to assume that $\beta_{11} = 0$. The automorphism from Example~6 with the scalar $\alpha = -\alpha_{12}$ allows us to assume that $\alpha_{12}=0$. Finally, the automorphism from Example~5 with the scalar $\alpha = \beta_{12}$ allows us to assume that $\beta_{12} = 0$. Thus, we can assume that $R(e_{22}) = -e_{22}$ or $R(e_{22}) = e_{11}$.

Let us introduce the following notation:
\begin{gather*}
 R(e_{21}) = \sum\gamma_{ij}e_{ij} + \sum\delta_{ij}ve_{ij},\\
 R(ve_{21}) = \sum\mu_{ij}e_{ij} + \sum\nu_{ij}ve_{ij},\\
 R(ve_{22}) = \sum\varepsilon_{ij}e_{ij} + \sum\eta_{ij}ve_{ij}.
\end{gather*}
Then
\[0 = R(e_{21})R(e_{11}) = R(R(e_{21})e_{11} + e_{21}) = R(\gamma_{21}e_{21} + e_{21}) = (\gamma_{21}+1)R(e_{21}),\]
whence $\gamma_{21}=-1$. Next,
\[0 = R(e_{11})R(e_{21}) = R(e_{11}R(e_{21})) = \delta_{21}R(ve_{21}) + \delta_{ 22}R(ve_{22}),\]
whence $\delta_{21}=\delta_{22}=0$. Next,
\[0 = R(ve_{21})R(e_{11}) = R(R(ve_{21})e_{11}) = \mu_{21}R(e_{21}),\]
whence $\mu_{21} = 0$. Next,
\[0 = R(e_{11})R(ve_{21}) = R(e_{11}R(ve_{21}) + ve_{21}) = (\nu_{21}+1)R( ve_{21}) + \nu_{22}R(ve_{22}),\]
whence $\nu_{22} = 0$, $\nu_{21} = -1$. Next,
\[0 = R(ve_{22})R(e_{11}) = R(R(ve_{22})e_{11}) = \varepsilon_{21}R(e_{21}),\]
whence $\varepsilon_{21} = 0$. Next,
\[0 = R(e_{11})R(ve_{22}) = R(e_{11}R(ve_{22}) + ve_{22}) = \eta_{21}R(ve_{21} ) + \eta_{22}R(ve_{22}) + R(ve_{22}),\]
whence $\eta_{21} = 0$, $\eta_{22} = -1$.

Let's consider two cases.

1) Let $R(e_{22}) = -e_{22}$. Then
\[ -\gamma_{12}e_{12}-\gamma_{22}e_{22} = R(e_{21})R(e_{22}) = R(R(e_{21})e_{22}) = \gamma_{22}R(e_{22}) = -\gamma_{22}e_{22},\]
whence $\gamma_{12} = 0$. Next,
\[e_{21}-\gamma_{22}e_{22}-\delta_{11}ve_{11}-\delta_{12}ve_{12} = R(e_{22})R(e_{21}) = R(-e_{21} + \gamma_{22}e_{22}) = -R(e_{21}) - \gamma_{22}e_{22},\]
whence $R(e_{21}) = -e_{21} + \delta_{11}ve_{11} + \delta_{12}ve_{12} \in -e_{21} + \mathrm{Ker}(R)$. Next,
\[-\mu_{22}e_{22}+ve_{21} = R(ve_{21})R(e_{22}) = R(R(ve_{21})e_{22}-ve_{21} + ve_{21}) = \mu_{22}R(e_{22})-R(ve_{21}),\]
whence $R(ve_{21}) = -ve_{21}$. Further,
\[-\varepsilon_{12}e_{12}-\varepsilon_{22}e_{22}- ve_{22} = R(ve_{22})R(e_{22}) = R(R(ve_{22})e_{22}) = -\varepsilon_{22}e_{22} - R(ve_{22}),\]
whence $R(ve_{22}) = -ve_{22} + \varepsilon_{12}e_{12} \in -ve_{22} + \mathrm{Ker}(R)$.

Thus, it is easy to see that $R(R(x)) = -R(x)$ for any $x\in\mathbb{O}$. Thus, the operator $R$ is splitting, it is a contradiction.

2) Let $R(e_{22}) = e_{11}$. Then
\[\gamma_{11}e_{11}+\gamma_{12}e_{12} = R(e_{22})R(e_{21}) = R(e_{22}R(e_{21}) + e_{21}) = \gamma_{22}R(e_{22}),\]
whence $\gamma_{12} = 0$, $\gamma_{11} = \gamma_{22}$. Further,
\[\gamma_{11}e_{11}-e_{21}+\delta_{11}ve_{11}+\delta_{12}ve_{12} = R(e_{21})R(e_{22}) = \gamma_{22}R(e_{22})+R(e_{21}),\]
whence $R(e_{21}) = -e_{21} + \delta_{11}ve_{11} + \delta_{12}ve_{12}$. The automorphism from Example~3 with the scalar $\alpha = \delta_{12}$ allows us to assume that $\delta_{12} = 0$, i.e. $R(e_{21}) = -e_{21} + \delta_{11}ve_{11}$. Next,
\[\mu_{11}e_{11}+\mu_{12}e_{12}-ve_{21}=R(e_{22})R(ve_{21})=R(ve_{21}+\mu_{22}e_{22}) = \mu_{22}e_{11} + R(ve_{21}),\]
whence $R(ve_{21}) = -ve_{21}+\mu_{11}e_{11}+\mu_{12}e_{12}$, $\mu_{22}=0$. Next,
\[\mu_{11}e_{11} = R(ve_{21})R(e_{22}) = R(\mu_{12}e_{12}) = 0,\]
whence $\mu_{11} = 0$ and $R(ve_{21}) = -ve_{21} + \mu_{12}e_{12}$. Next,
\[\varepsilon_{11}e_{11} + \varepsilon_{12}e_{12} - ve_{22} = R(e_{22})R(ve_{22}) = R(ve_{22} + \varepsilon_{22}e_{22}) = \varepsilon_{22}e_{11} + R(ve_{22}),\]
whence $R(ve_{22}) = -ve_{22} + \varepsilon_{11}e_{11} + \varepsilon_{12}e_{12}$ and $\varepsilon_{22}=0$. Next,
\[\varepsilon_{11}e_{11} = R(ve_{22})R(e_{22}) = R(-ve_{22} + \varepsilon_{12}e_{12} + ve_{22}) = 0,\]
whence $\varepsilon_{11} = 0$ and $R(ve_{22}) = -ve_{22} + \varepsilon_{12}e_{12}$. Further,
\[e_{21} - \mu_{12}ve_{12} + \varepsilon_{12}ve_{11} = R(ve_{22})R(ve_{21}) = R(-e_{21} - \varepsilon_{12}ve_{11} - e_{21} + \mu_{12}ve_{12} + e_{21}) = -R(e_{21}),\]
whence $\mu_{12} = 0$ and $\varepsilon_{12} = -\delta_{11}$. Example~4 with the scalar $\alpha = \delta_{11}$ allows us to assume that the operator $R$ coincides with $R_1$. The lemma is proved.

\section{Non-splitting Rota-Baxter operators of non-zero weight with three-dimensional kernel}

Let us prove that there are no non-splitting operators with idempotent three-dimensional kernel.

\medskip\textbf{Lemma 3.} \textit{Let $R$ be the Rota-Baxter operator of unit weight on the split Cayley-Dickson algebra $\mathbb{O}$ over $F$ and $\mathrm{Ker}(R)=Fe_{11}+Fve_{12}+Fve_{22}$. Then $R$ is splitting.
}

\textbf{Proof.} Let $R$ be the Rota-Baxter operator of unit weight on $\mathbb{O}$ with kernel $\mathrm{Ker}(R)=Fe_{11}+Fve_{12}+Fve_{22}$. Note that the theorem on the sum of the dimensions of the kernel and the image implies that the elements \[R(e_{12}),R(e_{21}),R(e_{22}),R(ve_{11}),R(ve_{21})\] form a basis of $\mathrm{Im}(R)$.

Let $R(e_{22}) = \sum\alpha_{ij}e_{ij} + \sum\beta_{ij}ve_{ij}$. Then
\[0=R(e_{22})R(e_{11})=R(\alpha_{11}e_{11}+\alpha_{21}e_{21}+\beta_{11}ve_{11}+\beta_{12}ve_{12})=\alpha_{21}R(e_{21})+\beta_{11}R(ve_{11}),\]
whence $\alpha_{21}=\beta_{11}=0$ due to the linear independence of the elements $R(e_{21})$ and $R(ve_{11})$. Further,
\[0 = R(e_{11})R(e_{22}) = R(\alpha_{11}e_{11}+\alpha_{12}e_{12}+\beta_{21}ve_{21}+\beta_{22}ve_{22}) = \alpha_{12}R(e_{12}) + \beta_{21}R(ve_{21}),\]
whence $\alpha_{12}=\beta_{21}=0$ due to the linear independence of the elements $R(e_{12})$ and $R(ve_{11})$. That is, $R(e_{22})=\alpha_{11}e_{11}+\alpha_{22}e_{22}+\beta_{12}ve_{12}+\beta_{22}ve_{22}$. Next,
\[\alpha_{11}^2e_{11} + (\alpha_{11}+\alpha_{22})\beta_{22}ve_{22} + \alpha_{22}^2e_{22} + \beta_{12}(\alpha_{11} + \alpha_{22})ve_{12} = R(e_{22})R(e_{22}) = (2\alpha_{22}+1)R(e_{22}),\]
where
\begin{gather*}
 \alpha_{11}(2\alpha_{22}-\alpha_{11}+1)=0,\quad \beta_{22}(\alpha_{22}-\alpha_{11}+1)=0,\quad \beta_{12}(\alpha_{22}-\alpha_{11}+1)=0,\quad (\alpha_{22}+1)\alpha_{22}=0.
\end{gather*}

Let us introduce the following notations:
\begin{gather*}
    R(e_{12}) = \sum\sigma_{ij}e_{ij} + \sum\tau_{ij}ve_{ij},\quad R(e_{21}) = \sum\gamma_{ij}e_{ij} + \sum\delta_{ij}ve_{ij},\\   R(ve_{21}) = \sum\mu_{ij}e_{ij} + \sum\nu_{ij}ve_{ij},\quad
    R(ve_{11}) = \sum\varepsilon_{ij}e_{ij} + \sum\eta_{ij}ve_{ij}.
\end{gather*}
Then
\[0 = R(e_{21})R(e_{11}) = R(R(e_{21})e_{11}+e_{21}) = R((\gamma_{21}+1)e_{21} + \delta_{11}ve_{11}) = (\gamma_{21}+1)R(e_{21})+\delta_{11}R(ve_{11}),\]
whence $\gamma_{21}=-1$, $\delta_{11}=0$ due to the linear independence of the elements $R(e_{21})$ and $R(ve_{11})$. Next,
\[0 = R(e_{11})R(e_{21}) = R(\gamma_{12}e_{12} + \delta_{21}ve_{21} + \delta_{22}ve_{22}) = \gamma_{12}R(e_{12}) + \delta_{21}R(ve_{21}),\]
whence $\gamma_{12} = \delta_{21} = 0$ due to the linear independence of the elements $R(e_{12})$ and $R(ve_{21})$. Next,
\[0 = R(e_{12})R(e_{11}) = R(R(e_{12})e_{11}) = R(\sigma_{21}e_{21}+\tau_{11}ve_{11}) = \sigma_{21}R(e_{21})+\tau_{11}R(ve_{11}),\]
whence $\sigma_{21}=\tau_{11}=0$ due to the linear independence of the elements $R(e_{21})$ and $R(ve_{11})$. Next,
\[0 = R(e_{11})R(e_{12}) = R(e_{11}R(e_{12})+e_{12}) = R((\sigma_{12}+1)e_{12}+\tau_{21}ve_{21}) = (\sigma_{12}+1)R(e_{12})+\tau_{21}R(ve_{21}),\]
whence $\sigma_{12}=-1$, $\tau_{21}=0$ due to the linear independence of the elements $R(e_{12})$ and $R(ve_{21})$. Next,
\[0 = R(ve_{21})R(e_{11}) = R(R(ve_{21})e_{11}) = R(\mu_{21}e_{21}+\nu_{11}ve_{11}) = \mu_{21}R(e_{21})+\nu_{11}R(ve_{11}),\]
whence $\mu_{21}=\nu_{11}=0$ due to the linear independence of the elements $R(e_{21})$ and $R(ve_{11})$. Next,
\[0 = R(e_{11})R(ve_{21}) = R(e_{11}R(ve_{21})+ve_{21}) = R(\mu_{12}e_{12} + (\nu_{21}+1)ve_{21}) = \mu_{12}R(e_{12}) + (\nu_{21}+1)R(ve_{21}),\]
whence $\mu_{12}=0$, $\nu_{21}=-1$ due to the linear independence of the elements $R(e_{12})$ and $R(ve_{21})$. Next,
\[0 = R(ve_{11})R(e_{11}) = R(R(ve_{11})e_{11} + ve_{11}) = R(\varepsilon_{21}e_{21} + (\eta_{11}+1)ve_{11}) = \varepsilon_{21}R(e_{21}) + (\eta_{11}+1)R(ve_{11}),\]
whence $\varepsilon_{21}=0$, $\eta_{11}=-1$ due to the linear independence of the elements $R(e_{21})$ and $R(ve_{11})$. Finally,
\[0 = R(e_{11})R(ve_{11}) = R(e_{11}R(ve_{11})) = R(\varepsilon_{12}e_{12} + \eta_{21}ve_{21}) = \varepsilon_{12}R(e_{12}) + \eta_{21}R(ve_{21}),\]
whence $\varepsilon_{12} = \eta_{21} = 0$ due to the linear independence of the elements $R(e_{12})$ and $R(ve_{21})$.

There are two possible cases.

1) $\alpha_{22}=0$. If $\alpha_{11}=0$, then $\beta_{22}=\beta_{12}=0$ and $R(e_{22})=0$, it is a contradiction. Thus, $\alpha_{11}\neq 0$, so that $\alpha_{11}=1$ and $R(e_{22})=e_{11}+\beta_{12}ve_{12}+\beta_{22}ve_{22}$. Example~2 with the scalar $\alpha = -\beta_{22}$ allows us to assume that $\beta_{22}=0$. Example~5 with the scalar $\alpha = \beta_{12}$ allows us to assume that $\beta_{12}=0$. Thus, we can assume that $R(e_{22})=e_{11}$.

Then 
\[\sigma_{11}e_{11} + \tau_{12}ve_{12} = R(e_{12})R(e_{22}) = R(-e_{12}+\sigma_{22}e_{22} + \tau_{22}ve_{22} + e_{12}) = \sigma_{22}e_{11},\]
whence $\tau_{12}=0$ and $\sigma_{11}=\sigma_{22}$. Next,
\[\sigma_{11}e_{11}-e_{12}+\tau_{22}ve_{22} = R(e_{22})R(e_{12}) = R(e_{12} + \sigma_{22}e_{22} + \tau_{12}ve_{12}) = R(e_{12}) + \sigma_{11}e_{11},\]
whence $R(e_{12})=-e_{12} + \tau_{22}ve_{22}$. Next,
\[\gamma_{11}e_{11} + \delta_{22}ve_{22} = R(e_{22})R(e_{21}) = R(-e_{21}+\gamma_{22}e_{22} + \delta_{12}ve_{12} + e_{21}) = \gamma_{22}e_{11},\]
whence $\gamma_{11}=\gamma_{22}$ and $\delta_{22}=0$. Next,
\[\gamma_{11}e_{11}-e_{21}+\delta_{12}ve_{12} = R(e_{21})R(e_{22}) = R(\gamma_{22}e_{22} + e_{21}) = \gamma_{11}e_{11} + R(e_{21}),\]
whence $R(e_{21}) = -e_{21}+\delta_{12}ve_{12}$. 
Next,
\[e_{11} = R(e_{12})R(e_{21}) = R(-e_{11} - e_{11} + e_{11}) = -R(e_{11}) = 0,\]
it is a contradiction. Thus, this case is impossible and $\alpha_{22}\neq 0$.

\medskip 2) $\alpha_{22}=-1$. If $\alpha_{11}=-1$, then $\beta_{22}=\beta_{12}=0$ and $R(e_{22})=-1$. By Proposition~1, the operator $R$ is splitting.

Let $\alpha_{11}\neq -1$, then $\alpha_{11}=0$. We have $R(e_{22})=-e_{22}+\beta_{12}ve_{12}+\beta_{22}ve_{22}$. Example~2 with the scalar $\alpha = -\beta_{22}$ allows us to assume that $\beta_{22}=0$. Example~5 with the scalar $\alpha = \beta_{12}$ allows us to assume that $\beta_{12}=0$. Thus, we can assume that $R(e_{22})=-e_{22}$.

Then
\[e_{21}-\gamma_{22}e_{22}-\delta_{12}ve_{12} = R(e_{22})R(e_{21}) = R(-e_{21}+\gamma_{22}e_{22} + \delta_{12}ve_{12} - e_{21} + e_{21}) = -\gamma_{22}e_{22} - R(e_{21}),\]
whence $R(e_{21}) = -e_{21} + \delta_{12}ve_{12}\in -e_{21}+\mathrm{Ker}(R)$. Next,
\[e_{12} - \sigma_{22}e_{22} = R(e_{12})R(e_{22}) = R(-e_{12}+\sigma_{22}e_{22}+\tau_{21}ve_{21} - e_{12} + e_{12}) = -R(e_{12}) - \sigma_{22}e_{22},\]
whence $R(e_{12})=-e_{12}$. Next,
\[-\varepsilon_{22}e_{22}+ve_{11}-\eta_{12}ve_{12} = R(e_{22})R(ve_{11}) = R(\varepsilon_{22}e_{22}-ve_{11}+\eta_{12}ve_{12}) = -\varepsilon_{22}e_{22}-R(ve_{11}),\]
whence $R(ve_{11})=-ve_{11}+\eta_{12}ve_{12}\in -ve_{11}+\mathrm{Ker}(R)$. Finally,
\[-\mu_{22}e_{22}+ve_{21}-\nu_{22}ve_{22} = R(ve_{21})R(e_{22}) = R(\mu_{22}e_{22}-ve_{21}+\nu_{22}ve_{22}) =
-\mu_{22}e_{22}-R(ve_{21}),\]
whence $R(ve_{21})=-ve_{21}+\nu_{22}ve_{22}\in -ve_{21}+\mathrm{Ker}(R)$.

Thus, it is easy to see that $R(R(x))=-R(x)$ for any $x\in\mathbb{O}$. Thus, the operator $R$ is splitting. The lemma is proved.

\medskip Let us construct an example of a non-splitting Rota-Baxter operator on $\mathbb{O}$ with a three-dimensional kernel.

\medskip\textbf{Example 20.} \textit{Let $\mathbb{O}$ be a split Cayley-Dickson algebra and $R_2:\mathbb{O}\to\mathbb{O}$ be a linear map defined as follows:
\[
    R_2(e_{21}) = -e_{21}, R_2(e_{22}) = -e_{22}, R_2(ve_{11}) = -ve_{11},  R_2(ve_{21}) =  -ve_{21},  R_2(e_{11}) = -1, R_2(Fe_{12}+Fve_{12}+Fve_{22})=0.
\]
Then $R_2$ is the Rota-Baxter operator on an algebra $\mathbb{O}$.
}

\smallskip It turns out that the operator $R_2$ is the only non-splitting operator on $\mathbb{O}$ with a three-dimensional kernel.

\medskip\textbf{Lemma 4.} \textit{Let $R$ be a non-splitting Rota-Baxter operator of unit weight on the split Cayley-Dickson algebra $\mathbb{O}$ over $F$ and $\mathrm{Ker}(R)=Fe_{12}+Fve_{12}+Fve_{22}$. Then, up to automorphism conjugation, the operator $R$ coincides with $R_2$.
}

\textbf{Proof.} Let $R$ be the Rota-Baxter operator of unit weight on $\mathbb{O}$ with kernel $\mathrm{Ker}(R)=Fe_{12}+Fve_{12}+Fve_{22}$. Note that the theorem on the sum of the dimensions of the kernel and the image implies that the elements \[R(e_{11}),R(e_{21}),R(e_{22}),R(ve_{11}),R(ve_{21})\] form a basis of $\mathrm{Im}(R)$.

Let $U = Fe_{11}+Fe_{21} + Fe_{22} + Fve_{11} + Fve_{21}$. Let $R(x)=1$, and we can assume that $x\in U$. Then
\[1 = R(x)R(x) = R(R(x)x+xR(x) + x^2) = 2R(x) + R(x^2) = 2 + R(x^2),\]
whence $R(x^2) = -1$. This means that $x+x^2\in\mathrm{Ker}(R)$. However, $U$ is a subalgebra and $U\cap\mathrm{Ker}(R) = 0$, so $x+x^2 = 0$. Let $x = \alpha_{11}e_{11} + \alpha_{21}e_{21} + \alpha_{22}e_{22} + \beta_{11}ve_{11} + \beta_{21}ve_{21}$. Then the condition $x^2 = -x$ implies
\begin{gather*}
\alpha_{11}^2 = -\alpha_{11},\\
\alpha_{21}(\alpha_{11}+\alpha_{22})=-\alpha_{21},\\
\alpha_{22}^2 = -\alpha_{22},\\
\beta_{11}(\alpha_{11}+\alpha_{22})=-\beta_{11},\\
\beta_{21}(\alpha_{11}+\alpha_{22})=-\beta_{21}.
\end{gather*}
If $\alpha_{11}=\alpha_{22}=0$, then $\alpha_{21}=\beta_{11}=\beta_{21}=0$ and $x=0$, it is a contradiction.

If $\alpha_{11}=\alpha_{22}=-1$, then $\alpha_{21}=\beta_{11}=\beta_{21}=0$ and $x=-1$, whence $R(1)=-1$ and the operator $R$ is splitting by Proposition~1.

Thus, either $\alpha_{11}=-1$, $\alpha_{22}=0$, or $\alpha_{11}=0$, $\alpha_{22}=-1$. This means that \[x = -e_{ii} + \alpha_{21}e_{21} + \beta_{11}ve_{11} + \beta_{21}ve_{21}\]
for some $i\in\{1,2\}$. But
\[R(y)=R(x)R(y)=R(y+xR(y)+xy) = R(y) + R(xR(y)) + R(xy),\]
whence $x(y+R(y))\in\mathrm{Ker}(R)$ for any $y\in\mathbb{O}$. In particular, $xe_{12}\in\mathrm{Ker}(R)$, i.e.
\[xe_{12}=-e_{ii}e_{12}+\alpha_{21}e_{22} + \beta_{21}ve_{11} \in Fe_{12}+Fve_{12}+Fve_{22},\]
whence $\alpha_{21}=\beta_{21}=0$. Further, $xve_{22}\in\mathrm{Ker}(R)$, i.e.
\[xve_{22} = -e_{ii}ve_{22} + \beta_{11}e_{22} \in Fe_{12}+Fve_{12}+Fve_{22},\]
whence $\beta_{11}=0$. Thus $x = -e_{ii}$ for some $i\in\{1,2\}$. Up to conjugation by canonical involution, we can assume that $x=-e_{11}$.

Thus, $R(e_{11})=-1$ and $e_{11}(y+R(y))\in\mathrm{Ker}(R)$ for any $y\in\mathbb{O}$. Similarly, $(y+R(y))e_{11}\in\mathrm{Ker}(R)$. Substituting $y\in\{e_{21},e_{22},ve_{11},ve_{21}\}$ into these inclusions , we obtain
\begin{gather*}
R(e_{21}) = -e_{21} + \delta_1 e_{22} + \varepsilon_1e_{12} + \varepsilon_2 ve_{12} + \varepsilon_3 ve_{22},\\
R(ve_{11}) = -ve_{11} + \delta_2 e_{22} + \eta_1 e_{12} + \eta_2 ve_{12} + \eta_3 ve_{22},\\
R(ve_{21}) = -ve_{21} + \delta_3 e_{22} + \sigma_1 e_{12} + \sigma_2 ve_{12} + \sigma_3 ve_{22},\\
 R(e_{22}) = \delta_4 e_{22} + \tau_1 e_{12} + \tau_2 ve_{12} + \tau_3 ve_{22}
\end{gather*}
for some $\delta_i, \varepsilon_i, \eta_i, \sigma_i, \tau_i\in F$. Then $n(R(e_{11}-e_{22})) = n(-1-\delta_4e_{22} - \tau_1 e_{12}-\tau_2ve_{12}-\tau_3 ve_{22}) = 1+\delta_4.$ By (\cite{BGP}, Lemma 3.1.c) either $R$ is a splitting operator or $1+\delta_4 = 0$, that is, $\delta_4 = -1$.

Next,
\[-\varepsilon_1\cdot 1 + \delta_1 R(e_{21}) = R(e_{21})R(e_{21}) = R(\varepsilon_1\cdot 1 + \delta_1 e_{21}) = \varepsilon_1 R(1) + \delta_1 R(e_{21}),\]
whence $\varepsilon_1 = 0$ (since $R(1)\neq -1$). Next,
\[- \eta_3\cdot 1 +\delta_2 R(ve_{11}) = R(ve_{11})R(ve_{11}) = R(\eta_3\cdot 1 + \delta_2 ve_{11} ) = \eta_3R(1) + \delta_2 R(ve_{11}),\]
whence $\eta_3 = 0$ (since $R(1)\neq -1$). Next,
\[-\sigma_2\cdot 1 + \delta_3 R(ve_{21}) = R(ve_{21})R(ve_{21}) = R(\sigma_2\cdot 1 + \delta_3 ve_{21} ) = \sigma_2 R(1) + \delta_3 R(ve_{21}),\]
whence $\sigma_2 = 0$ (since $R(1)\neq -1$). Next,
\[-\tau_1e_{11} + \tau_1\delta_1 e_{12} - e_{21} - \delta_1e_{22} - (\varepsilon_2 + \tau_1\varepsilon_3)ve_{12} + (\tau_3\delta_1 - \tau_2)ve_{22} = R(e_{22})R(e_{21}) = -R(e_{21}) + \tau_1 R(e_{11}) + \delta_1 R(e_{22 }),\]
from which we obtain (comparing the coefficients of $e_{22}$ and $ve_{22}$) $\tau_1=-\delta_1$, $\tau_2=\varepsilon_3$.
Next,
\begin{multline*}
-\eta_1 e_{11} + (\eta_1\delta_1-\varepsilon_2)e_{12} -\delta_2 e_{21} + (\delta_1\delta_2 - \varepsilon_3)e_{22} + (\delta_2\varepsilon_2 - \ eta_1\varepsilon_3)ve_{12} + ve_{21} - \eta_2ve_{22}
=\\ = R(ve_{11})R(e_{21}) = -R(ve_{21}) + \delta_2 R(e_{21}) + \eta_1 R(e_{11}) + \varepsilon_3 R(e_{22}),
\end{multline*}
so we get (comparing coefficients of $e_{22}$ and $ve_{22}$) $\delta_3 = -\eta_1$ and
\begin{equation}\label{etasigma}
-\eta_2 = -\sigma_3 + \delta_2\varepsilon_3 + \varepsilon_3\tau_3.
\end{equation}
Next,
\begin{multline*}
-\varepsilon_3 e_{11} + \varepsilon_2 e_{12} + (\delta_1\delta_2 - \eta_1) e_{22} - \delta_1 ve_{11} + (\delta_1\eta_2 + \varepsilon_3\eta_1)ve_{12 } - ve_{21} + (\eta_2 + \varepsilon_3\delta_2)ve_{22} = \\ = R(e_{21})R(ve_{11}) = \varepsilon_3 R(e_{11}) + \eta_1 R(e_{22}) + \delta_1 R(ve_{11}) + R(ve_{21}).
\end{multline*}
It means that (we compare coefficients of $e_{12}$ and $e_{22}$) $\delta_3 = \varepsilon_3$ and $\varepsilon_2 = \sigma_1$. Next,
\[-\tau_3 e_{11} + (\delta_3 - \delta_1\delta_2)e_{12} - \delta_2 e_{22} + ve_{11} - (\eta_2+\tau_3\delta_3)ve_{12} + \tau_3\delta_2 ve_{22} = R(e_{22})R(ve_{11}) = -R(ve_{11}) + \tau_3 R(e_{11}) + \delta_2 R(e_{ 22}).\]
It means that (we compare coefficients of $e_{12}$ and $e_{22}$) $\delta_2 = -\delta_3 = -\tau_3$. Then from \eqref{etasigma} we obtain $\eta_2 = \sigma_3$.

Example~3 with scalar $\alpha = \sigma_1$ allows us to assume that $\sigma_1 = 0$. Example~1 with scalar $\alpha = \sigma_3$ allows us to assume that $\sigma_3 = 0$. Thus,
\begin{gather*}
 R(e_{21}) = -e_{21} + \delta_1 e_{22} - \delta_2 ve_{22},\quad
 R(ve_{11}) = -ve_{11} + \delta_2 e_{22} + \delta_2 e_{12} ,\\
 R(ve_{21}) = -ve_{21} - \delta_2 e_{22},\quad
 R(e_{22}) = - e_{22} - \delta_1 e_{12} - \delta_2 ve_{12} - \delta_2 ve_{22}.
\end{gather*}

Suppose that $\delta_2 \neq 0$. Then Example~7 with the scalar $\alpha = \delta_2$ allows us to assume that $\delta_2 = 1$, i.e.
\begin{gather*}
R(e_{21}) = -e_{21} + \delta_1 e_{22} - ve_{22},\quad
R(ve_{11}) = -ve_{11} + e_{22} + e_{12},\\
R(ve_{21}) = -ve_{21} - e_{22},\quad
R(e_{22}) = - e_{22} - \delta_1 e_{12} - ve_{12} - ve_{22}.
\end{gather*}
However, Example~6 with the scalar $\alpha = 1$ allows us to assume that $R(ve_{11})=-ve_{11}+e_{22}$. This contradicts the above proof that the coefficients in $R(ve_{11})$ for $e_{22}$ and $e_{12}$ must coincide. Thus, $\delta_2 = 0$, that is
\begin{gather*}
R(e_{21}) = -e_{21} + \delta_1 e_{22},\quad R(ve_{11}) = -ve_{11}, \quad R(ve_{21}) = -ve_{21},\quad
R(e_{22}) = - e_{22} - \delta_1 e_{12}.
\end{gather*}
Example~6 with scalar $\alpha = \delta_1$ allows us to assume that $\delta_1 = 0$, that is,
\[
R(e_{21}) = -e_{21},\quad R(ve_{11}) = -ve_{11},\quad R(ve_{21}) = -ve_{21},\quad R(e_{22}) = - e_{22}.
\]
The lemma is proved.

\medskip Now we are ready to obtain a complete description of non-splitting operators on $\mathbb{O}$.

\medskip\textbf{Definition.} We introduce the following notation: if $R$ is a Rota-Baxter operator of weight $\lambda$ on an algebra $A$, then $\phi(R) = -R - \lambda\mathrm{id}$, where $\phi(R):A\to A$.

It is easy to see that the transformation $\phi(R)$ is also a Rota-Baxter operator of weight $\lambda$ on an algebra $A$.

\medskip Note that the subalgebra $B_1 = Fe_{11} + Fe_{12} + Fve_{11} + Fve_{12}$ can be translated by an automorphism of the algebra $\mathbb{O}$ into the subalgebra $Fe_{11} + Fe_{12} + Fve_{12} + Fve_{22}$, and the subalgebra $Fe_{11}+Fe_{22}$ remains fixed by this automorphism. In this case, the operator $R_1$ goes over to the operator $\phi(R_2)$. Denote $B_2 = Fe_{21} + Fve_{21} + Fve_{22}$.

\medskip\textbf{Theorem 1.} \textit{Let $R$ be a non-splitting Rota-Baxter operator of non-zero weight on a split Cayley-Dickson algebra $\mathbb{O}$. Then, up to conjugation by an automorphism, an anti-automorphism, up to multiplication by a non-zero scalar, and up to the action of $\phi$, the operator $R$ coincides with the operator~$R_1$:
\[R_1|_{B_1} = -\mathrm{id}, \quad R_1|_{B_2} = 0, \quad R_1(e_{22}) = e_{11}.\]
}

\textbf{Proof.} By Lemma~1, $\dim(\mathrm{Ker}(R)) \in\{3,4\}$. Up to multiplication by a scalar, we may assume that the operator $R$ has unit weight. By Proposition~1, we may assume that $\mathrm{Ker}(R)$ does not contain the identity of $\mathbb{O}$. By the results of \cite{Octonions}, up to automorphism and antiautomorphism, there exist exactly two non-unital three-dimensional subalgebras: $Fe_{11}+Fve_{12}+Fve_{22}$ and $Fe_{12}+Fve_{12}+Fve_{22}$. Then Lemma~2, Lemma~3, and Lemma~4 imply that $R=R_1$ or $R=R_2$ up to automorphism and antiautomorphism. But $\phi(R_2)=R_1$. The theorem is proved.

\medskip\textbf{Definition.} A unital algebra $A$ over a field $F$ with a quadratic form $n:A\to F$ is called a \textbf{\textit{composition algebra}} if $n(xy)=n(x)n(y)$ for any $x,y\in A$ and a bilinear form $f(x,y)=n(x+y)-n(x)-n(y)$ is nondegenerate. 

A composition algebra with zero divisors is called \textbf{\textit{split}}. It is well known that there are only 3 split composition algebras: $F\oplus F$, $M_2(F)$, $\mathbb{O}$. The papers \cite{BGP} and \cite{Unital} allow us to state the following

\medskip\textbf{Corollary 1.} {\it In any split composition algebra there is the unique Rota-Baxter operator of non-zero weight up to conjugation by an automorphism, an anti-automorphism, up to multiplication by a non-zero scalar, and up to the action of~$\phi$.}

    \section{Decomposition of split octonions into a sum of two subalgebras}

\textbf{Lemma 5.} {\it Let $\mathbb{O} = A_1 \oplus A_2$ be a decomposition of the algebra of split octonions into a direct sum of two subalgebras and $\dim A_1 = 6$. Then, up to automorphism and antiautomorphism, we can assume that $A_1 = M_2(F)+Fve_{12}+Fve_{22}$, and one of the following cases holds for the algebra $A_2$:
\begin{enumerate}
\item $A_2 = Fve_{11} + Fve_{21}$;
\item $A_2 = Fve_{11} + F(ve_{21}+e_{22})$.
\end{enumerate}
}

\textbf{Proof.} Let $\mathbb{O} = A_1 \oplus A_2$, where $\dim A_1 = 6$. By the results of \cite{Octonions}, we can assume that $A_1 = M_2(F) + Fve_{12} + Fve_{22}$. Then the basis in $A_2$ is of the form
\[x_1 = ve_{11} + \sum\alpha_{ij}e_{ij} + \gamma_1 ve_{12} + \gamma_2 ve_{22}, \quad x_2 = ve_{21} + \sum\beta_{ij}e_{ij} + \delta_1 ve_{12} + \delta_2 ve_{22}.\]
Then $x_1^2\in A_2$. But the coefficient of $ve_{11}$ in $x_1^2$ is $\alpha_{11} + \alpha_{22}$, and the coefficient of $ve_{21}$ is $0$, so $x_1^2 = (\alpha_{11}+\alpha_{22})x_1$. Similarly, $x_2^2 = (\beta_{11} + \beta_{22}) x_2$.

The element $x_1^2$ has the coefficient of $e_{11}$ to be $\alpha_{11}^2+\alpha_{12}\alpha_{21} + \gamma_2$. Using the equality $(\alpha_{11}+\alpha_{22})x_1 = x_1^2$, we obtain
\[\gamma_2 = \alpha_{11}\alpha_{22}-\alpha_{12}\alpha_{21}.\]
Similarly, we obtain
\[\delta_1 = -\beta_{11}\beta_{22}+\beta_{12}\beta_{21}.\]
The element $x_1x_2$ has the coefficient of $ve_{11}$ to be $\beta_{11}-\alpha_{12}$, and the coefficient of $ve_{21}$ is equal to $\beta_{21}+\alpha_{11}$. This means that $x_1x_2 = (\beta_{11}-\alpha_{12})x_1 + (\beta_{21}+\alpha_{11})x_2$. For $x_2x_1$, the coefficient of $ve_{11}$ is $\alpha_{12}+\beta_{22}$, and the coefficient of $ve_{21}$ is $\alpha_{22} - \beta_{21}$. This means that $x_2x_1 = (\alpha_{12}+\beta_{22})x_1 + (\alpha_{22}-\beta_{21})x_2$. Note that by the results of \cite{Octonions} we can assume (up to automorphism or antiautomorphism) that $A_2\simeq Fe_{11}+Fe_{12}$ or $A_2$ has zero multiplication. In both cases $ab\in Fb$ for any $a,b\in A_2$. Then from the equalities above it follows that $\alpha_{22}=\beta_{21}$, $\beta_{11}=\alpha_{12}$, $x_1x_2 = (\alpha_{11}+\alpha_{22})x_2$, $x_2x_1 = (\alpha_{12} + \beta_{22})x_1$. From the last condition, presented on the basis vectors, we obtain $-\alpha_{11}\beta_{22}+\alpha_{21}\beta_{12}-\gamma_1+\delta_2 = 0$.
Example~3 with the scalar $\alpha = -\beta_{12}$ allows us to assume that
\begin{gather*}
x_1 = ve_{11} + \alpha_{11}e_{11} + \alpha_{12}e_{12} + \alpha_{21}e_{21} + \alpha_{22}e_{22} + \gamma_3 ve_{12} + (\alpha_{11}\alpha_{22} - \alpha_{12}\alpha_{21})ve_{22},\\
x_2 = ve_{21} + \alpha_{12}e_{11} + \alpha_{22}e_{21} + \beta_{22}e_{22} -\alpha_{12}\beta_{22}ve_{12} + (\gamma_3 + \alpha_{11}\beta_{22}) ve_{22}.
\end{gather*}
where $\gamma_3 = \gamma_1 - \beta_{12}\alpha_{21}.$ Example~1 with a scalar $\alpha = -\gamma_3$ allows us to assume that $\gamma_3 = 0$, that is
\begin{gather*}
x_1 = ve_{11} + \alpha_{11}e_{11} + \alpha_{12}e_{12} + \alpha_{21}e_{21} + \alpha_{22}e_{22} + (\alpha_{11}\alpha_{22} - \alpha_{12}\alpha_{21})ve_{22},\\
x_2 = ve_{21} + \alpha_{12}e_{11} + \alpha_{22}e_{21} + \beta_{22}e_{22} -\alpha_{12}\beta_{22}ve_{12} + \alpha_{11}\beta_{22}ve_{22}.
\end{gather*}
Example~8 with the scalar $\alpha = -\alpha_{21}$ allows us to assume that
\begin{gather*}
x_1 = ve_{11} + \alpha_{11}e_{11} + \alpha_{12}e_{12} + \alpha_{22}e_{22} + \alpha_{11}\alpha_{22}ve_{22},\\
x_2 = ve_{21} + \alpha_{12}e_{11} + \alpha_{22}e_{21} + \beta_{22}e_{22} -\alpha_{12}\beta_{22}ve_{12} + \alpha_{11}\beta_{22} ve_{22}.
\end{gather*}
Example~2 with the scalar $\alpha = \alpha_{11}$ allows us to assume that
\begin{gather*}
x_1 = ve_{11} + \alpha_{12}e_{12} + (\alpha_{11}+\alpha_{22})e_{22} + \alpha_{11}\alpha_{12}ve_{12},\\
x_2 = ve_{21} + \alpha_{12}e_{11} + (\alpha_{11}+\alpha_{22})e_{21} + \beta_{22}e_{22} -\alpha_{12}\beta_{22}ve_{12} + \alpha_{11}\alpha_{12} ve_{22}.
\end{gather*}
Example~5 with the scalar $\alpha = \alpha_{12}$ allows us to assume that $x_1 = ve_{11} + (\alpha_{11}+\alpha_{22})e_{22} + \alpha_{12}(2\alpha_{11}+\alpha_{22})ve_{12}$, $x_2 = ve_{21} + (\alpha_{11}+\alpha_{22})e_{21} + (\alpha_{12}+\beta_{22})e_{22} + \alpha_{12}(2\alpha_{11}+\alpha_{22})ve_{22}$. Example~1 with the scalar $\alpha = -\alpha_{12}(2\alpha_{11}+\alpha_{22})$ allows us to assume that $x_1 = ve_{11} + (\alpha_{11}+\alpha_{22})e_{22}$, $x_2 = ve_{21} + (\alpha_{11}+\alpha_{22})e_{21} + (\alpha_{12}+\beta_{22})e_{22}$. Let's consider several cases.

1) Let $\alpha_{11}+\alpha_{22}=\alpha_{12}+\beta_{22} = 0$. Then $x_1 = ve_{11}$, $x_2 = ve_{21}$.

2) Let $\alpha_{11}+\alpha_{22} = 0$, $\alpha_{12} + \beta_{22} \neq 0$. Then $x_1 = ve_{11}$, $x_2 = ve_{21} + (\alpha_{12}+\beta_{22})e_{22}$. Example~7 with the scalar $\alpha = (\alpha_{12}+\beta_{22})$ and multiplication by $\frac{1}{\alpha}$ allows us to assume that $x_1 = ve_{11}$, $x_2 = ve_{21} + e_{22}$.

3) Let $\alpha_{12}+\beta_{22} = 0$, $\alpha_{11}+\alpha_{22}\neq 0$. Then $x_1 = ve_{11} + (\alpha_{11}+\alpha_{22})e_{22}$, $x_2 = ve_{21} + (\alpha_{11}+\alpha_{22})e_{21}$. Example~7 with the scalar $\alpha = -(\alpha_{11}+\alpha_{22})$ and multiplication by $-\frac{1}{\alpha}$ allow us to assume that $x_1 = ve_{11} - e_{22}$, $x_2 = ve_{21} - e_{21}$. Example~2 with the scalar $\alpha = 1$ allows us to assume that $x_1 = ve_{11} - e_{11}$, $x_2 = ve_{21}$. The antiautomorphism from Example~11 allows us to assume that $x_1 = -ve_{21} - e_{11}$, $x_2 = ve_{11}$. Classical involution allows us to assume that $x_1 = ve_{21} - e_{22}$, $x_2 = ve_{11}$. Example~7 with the scalar $\alpha = -1$ allows us to assume that $x_1 = ve_{21} + e_{22}$, $x_2 = ve_{11}$.

4) Let $\alpha_{11}+\alpha_{22}\neq 0$, $\alpha_{12}+\beta_{22}\neq 0$. Example~7 with the scalar $\alpha = -(\alpha_{11}+\alpha_{22})$ and multiplication by $-\frac{1}{\alpha}$ allow us to assume that $x_1 = ve_{11} - e_{22}$, $x_2 = ve_{21} - e_{21} + \delta e_{22}$, where $\delta = -\frac{\alpha_{12}+\beta_{22}}{\alpha_{11}+\alpha_{22}}\neq 0$. Example~11 allows us to assume that $x_1 = ve_{21}+e_{22}$, $x_2 = ve_{11}-e_{12}+\delta e_{22}$. Example~12 with the scalar $\alpha = \delta$ allows us to assume that $x_2 = ve_{11} - e_{12} + e_{22}$. Example~17 allows us to assume that $x_1 = ve_{11}$, $x_2 = ve_{21}+e_{22}$

The lemma is proved.

\medskip\textbf{Lemma 6.} {\it Let $\mathbb{O} = A_1 \oplus A_2$ be the decomposition of the algebra of split octonions into a direct sum of two subalgebras, where $\dim A_1 = 5$. Then, up to automorphism and antiautomorphism, we can assume that $A_1 = Fe_{11} + Fe_{12} + Fe_{22} + Fve_{12} + Fve_{22}$, and for the algebra $A_2$ one of the following options holds:
\begin{enumerate}
\item $A_2 = Fe_{21} + Fve_{11} + Fve_{21}$,
\item $A_2 = F(e_{21}+e_{11}) + Fve_{11} + Fve_{21}$,
\end{enumerate}
}

\textbf{Proof.} Let $\mathbb{O} = A_1 \oplus A_2$, where $\dim A_1 = 5$. By the results of \cite{Octonions}, we can assume that $A_1 = Fe_{11} + Fe_{12} + Fe_{22} + Fve_{12} + Fve_{22}$. Then the basis in $A_2$ has the form
\begin{gather*}
 x_1 = e_{21} + \alpha_{11}e_{11} + \alpha_{12}e_{12} + \alpha_{22}e_{22} + \delta_1 ve_{12} + \delta_2 ve_{22},\\
 x_2 = ve_{11} + \beta_{11}e_{11} + \beta_{12}e_{12} + \beta_{22}e_{22} + \mu_1 ve_{12} + \mu_2 ve_{22},\\
 x_3 = ve_{21} + \gamma_{11}e_{11} + \gamma_{12}e_{12} + \gamma_{22}e_{22} + \nu_1 ve_{12} + \nu_2 ve_{22}.
\end{gather*}
Further, regardless of the use of automorphisms and antiautomorphisms, we maintain in our notations the situation when $x_1\in e_{21}+A_1$, $x_2\in ve_{11}+A_1$, $x_3\in ve_{21}+A_1$. Note that in this case, for a fixed subalgebra $A_2$, the elements $x_1,x_2,x_3$ are uniquely determined.

Then $x_1^2\in A_2$. But the coefficient of $e_{21}$ in $x_1^2$ is $\alpha_{11} + \alpha_{22}$, and the coefficients of $ve_{11}$ and $ve_{21}$ are $0$, so $x_1^2 = (\alpha_{11}+\alpha_{22})x_1$. Comparing the coefficients of $e_{11}$ in this equality, we get $\alpha_{12} = \alpha_{11}\alpha_{22}$. Similarly, $x_2^2 = (\beta_{11}+\beta_{22})x_2$ and $\mu_2=\beta_{11}\beta_{22}$, $x_3^2 = (\gamma_{11}+\gamma_{22})x_3$ and $\nu_1 = -\gamma_{11}\gamma_{22}$.

Since $x_1x_2\in \beta_{11}e_{21} + \alpha_{22}ve_{11} - ve_{21} + A_1$, then $x_1x_2 = \beta_{11}x_1 + \alpha_{22}x_2 - x_3$. Comparing the coefficients of $e_{11},e_{12},e_{22},ve_{22}$ in this equality, we get
\begin{gather*}
\gamma_{11} = -\delta_2 + \alpha_{22}\beta_{11},\\
\gamma_{12} = \alpha_{11}\alpha_{22}\beta_{11} - \alpha_{11}\alpha_{22}\beta_{22} - \alpha_{11}\beta_{12}+\alpha_{22}\beta_{12} + \delta_1,\\
\gamma_{22} = \alpha_{22}\beta_{11} + \beta_{12},\\
\nu_2 = -\alpha_{11}\beta_{11}\beta_{22}+\alpha_{22}\beta_{11}\beta_{22}+\delta_2\beta_{11}-\delta_2\beta_{22}+\mu_1.
\end{gather*}
Next, $x_2x_1\in \beta_{22}e_{21} + \alpha_{11}ve_{11} + ve_{21}$, so $x_2x_1 = \beta_{22}x_1 + \alpha_{11}x_2 + x_3$. Comparing the coefficients of $e_{12}$, we obtain
\[\beta_{12} = \alpha_{11}\beta_{22}+\alpha_{22}\beta_{11}-\delta_2.\]
Example~3 with the scalar $\alpha = -\delta_1$ allows us to assume that $\delta_1 = 0$. Example~1 with the scalar $\alpha = -\mu_1$ allows us to assume that $\mu_1 = 0$. Thus, 
\begin{gather*}
    x_1 = e_{21} + \alpha_{11}e_{11} + \alpha_{11}\alpha_{22}e_{12} + \alpha_{22}e_{22} + \delta_2 ve_{22},\\
    x_2 = ve_{11} + \beta_{11}e_{11} + (\alpha_{11}\beta_{22}+\alpha_{22}\beta_{11}-\delta_2)e_{12} + \beta_{22}e_{22} + \beta_{11}\beta_{22} ve_{22},\\
    x_3 = ve_{21} - (\delta_2 - \alpha_{22}\beta_{11})e_{11} + (\alpha_{11}\delta_2 - \alpha_{22}\delta_2-\alpha_{11}^2\beta_{22}+\alpha_{22}^2\beta_{11})e_{12} + (\delta_2-\alpha_{11}\beta_{22})e_{22}  + \\ (\delta_2 - \alpha_{22}\beta_{11})(\delta_2 -\alpha_{11}\beta_{22}) ve_{12} + (-\alpha_{11}\beta_{11}\beta_{22}+\alpha_{22}\beta_{11}\beta_{22}+\delta_2\beta_{11}-\delta_2\beta_{22}) ve_{22}.
\end{gather*}

Let us consider several cases.

1) Let $\alpha_{11}\neq 0$. Example~13 with a scalar $\alpha = \frac{1}{\alpha_{11}}$ and multiplying elements $x_1,x_3$ by the scalar $\frac{1}{\alpha_{11}}$ allows us to assume that
\begin{gather*}
x_1 = e_{21} + e_{11} + \alpha_{22}'e_{12} + \alpha_{22}'e_{22} + \delta_2've_{22},\\
x_2 = ve_{11} + \beta_{11}e_{11} + (\beta_{22}+\alpha_{22}'\beta_{11}-\delta_2')e_{12} + \beta_{22}e_{22} + \beta_{11}\beta_{22}ve_{22},\\
 x_3 = ve_{21} - (\delta_2' - \alpha_{22}'\beta_{11})e_{11} + (\delta_2'- \alpha_{22}'\delta_2'-\beta_{22}+\alpha_{22}'^2\beta_{11})e_{12} + (\delta_2'-\beta_{22})e_{22} + \\ (\delta_2' - \alpha_{22}'\beta_{11})(\delta_2'-\beta_{22}) ve_{12} + (-\beta_{11}\beta_{22}+\alpha_{22}'\beta_{11}\beta_{22}+\delta_2'\beta_{11}-\delta_2'\beta_{22}) ve_{22},
\end{gather*}
where $\delta_2' = \frac{\delta_2}{\alpha_{11}}$, $\alpha_{22}'=\frac{\alpha_{22}}{\alpha_{11}}$. The automorphism from Example~6 with a scalar $\alpha = -\alpha_{22}'$ and considering $x_3-\alpha_{22}'x_2$ instead of $x_3$ allows us to assume that
\begin{gather*}
x_1 = e_{21} + (\alpha_{22}'+1)e_{11} +\alpha_{22}'\delta_2' ve_{12} + \delta_2' ve_{22},\\
x_2 = ve_{11} + \beta_{11}e_{11} + (\beta_{22}+\alpha_{22}'\beta_{22}-\delta_2')e_{12} + \beta_{22}e_{22} + \alpha_{22}'\beta_{11}\beta_{22}ve_{12} + \beta_{11}\beta_{22}ve_{22},\\
 x_3 = ve_{21} - \delta_2'e_{11} + (-\alpha_{22}'^2\beta_{22}+2\alpha_{22}'\delta_2'-2\alpha_{22}'\beta_{22}+\delta_2'-\beta_{22})e_{12} + (-\beta_{22}+\delta_2'-\alpha_{22}'\beta_{22})e_{22} +\\ \delta_2'(-\alpha_{22}'\beta_{22}+\delta_2'-\beta_{22}) ve_{12} + (\delta_2'\beta_{11}-\delta_2'\beta_{22}-\beta_{11}\beta_{22}) ve_{22}.
\end{gather*}

1.a) Let $\delta_2'=0$. Then
\begin{gather*}
 x_1 = e_{21} + (\alpha_{22}'+1)e_{11},\\
 x_2 = ve_{11} + \beta_{11}e_{11} + \beta_{22}(\alpha_{22}'+1)e_{12} + \beta_{22}e_{22} + \alpha_{22}'\beta_{11}\beta_{22}ve_{12} + \beta_{11}\beta_{22}ve_{22},\\
 x_3 = ve_{21} - \beta_{22}(\alpha_{22}'+1)^2e_{12} - \beta_{22}(\alpha_{22}'+1)e_{22} -\beta_{11}\beta_{22}ve_{22}.
\end{gather*}

1.a.a) Let $\beta_{11}\neq 0$. Then Example~7 with a scalar $\alpha = \beta_{11}$ and multiplying elements $x_2,x_3$ by scalar $\frac{1}{\beta_{11}}$ allow us to assume that
\begin{gather*}
x_1 = e_{21} + (\alpha_{22}'+1)e_{11},\\
x_2 = ve_{11} + e_{11} + \beta(\alpha_{22}'+1)e_{12} + \beta e_{22} + \alpha_{22}'\beta ve_{12} + \beta ve_{22},\\
x_3 = ve_{21} - \beta(\alpha_{22}'+1)^2e_{12} - \beta(\alpha_{22}'+1)e_{22} -\beta ve_{22},
\end{gather*}
where $\beta = \frac{\beta_{22}}{\beta_{11}}$.

1.a.a.a) Let $\alpha_{22}'+1\neq 0$. Then Example~8 with a scalar $\alpha = \frac{1}{\alpha_{22}'+1}$, consideration an element $x_2-\frac{1}{\alpha_{22}'+1}x_1$ instead of an element $x_2$, Example~1 with scalar $\alpha = -\alpha_{22}'\beta$, an automorphism from Example~13 with a scalar $\alpha =\frac{1}{\alpha_{22}'+1}$, multiplying $x_1,x_3$ by $\frac{1}{\alpha_{22}'+1}$ allow us to assume that
\begin{gather*}
x_1 = e_{21} + e_{11},\\
x_2 = ve_{11} + \beta e_{12} + \beta e_{22},\\ 
x_3 = ve_{21} -\beta e_{12} -\beta e_{22}.
\end{gather*}

1.a.a.a.a) Let $\beta = 0$. Then $x_1 = e_{21}+e_{11}$, $x_2 = ve_{11}$, $x_3 = ve_{21}$.

1.a.a.a.b) Let $\beta\neq 0$. Then the automorphism from Example~7 with a scalar $\alpha = \beta$ and the multiplication of $x_2,x_3$ by the scalar $\frac{1}{\beta}$ allow us to assume that
\begin{gather*}
x_1 = e_{21} + e_{11},\\
x_2 = ve_{11} + e_{12} + e_{22},\\
x_3 = ve_{21} - e_{12} - e_{22}.
\end{gather*}
Example~6 with a scalar $\alpha = 1$ and consideration an element $x_3+x_2$ instead of an element $x_3$ allow us to assume that
\begin{gather*}
x_1 = e_{21} + e_{22},\\
x_2 = ve_{11} + e_{22},\\
x_3 = ve_{21}.
\end{gather*}
Consideration an element $x_2 - x_1$ instead of an element $x_2$ and Example~8 with scalar $\alpha = 1$ give us $x_1 = e_{21}+e_{22}$, $x_2 = ve_{11}$, $x_3 = ve_{21}$. Then the canonical involution and Example~12 with the scalar $\alpha = -1$ allow us to assume that we are in case 1.a.a.a.a).

1.a.a.b) Let $\alpha_{22}'+1=0$. The automorphism from Example~2 with a scalar $\alpha =- \beta$ and the automorphism from Example~1 with a scalar $\alpha = \beta$ allow us to assume that
\begin{gather*}
x_1 = e_{21},\\
x_2 = ve_{11} + (\beta + 1) e_{11},\\
x_3 = ve_{21}.
\end{gather*}

1.a.a.b.a) Let $\beta = -1$. Then $x_1 = e_{21}$, $x_2 = ve_{11}$, $x_3 = ve_{21}$.

1.a.a.b.b) Let $\beta+1\neq 0$. Then Example~7 with a scalar $\alpha = \beta + 1$, multiplication of $x_2,x_3$ by scalar $\frac{1}{\beta + 1}$, Example~14, replacement $x_1\longleftrightarrow x_2$, canonical involution and Example~12 with a scalar $\alpha = -1$ allow us to consider that we are in case 1.a.a.a.a).

1.a.b) Let $\beta_{11} = 0$. Then
\begin{gather*}
 x_1 = e_{21} + (\alpha_{22}'+1)e_{11},\\
 x_2 = ve_{11} + \beta_{22}(\alpha_{22}'+1)e_{12} + \beta_{22}e_{22},\\
 x_3 = ve_{21} - \beta_{22}(\alpha_{22}'+1)^2e_{12} - \beta_{22}(\alpha_{22}'+1)e_{22}.
\end{gather*}

1.a.b.a) Let $\alpha_{22}'+1\neq 0$. Then Example~13 with a scalar $\alpha = \frac{1}{\alpha_{22}'+1}$ and the multiplication of elements $x_1,x_3$ by a scalar $\frac{1}{\alpha_{22}'+1}$ allow us to assume that we have fallen into case 1.a.a.a).

1.a.b.b) Let $\alpha_{22}'+1=0$.

1.a.b.b.a) Let $\beta_{22}\neq 0$. Then Example~7 with a scalar $\alpha = \beta_{22}$ and multiplication of elements $x_2,x_3$ by a scalar $\frac{1}{\beta_{22}}$ allow us to assume that $x_1 = e_{21}$, $x_2 = ve_{11} + e_{22}$, $x_3 = ve_{21}$. Similarly to 1.a.a.b.b), this case reduces to 1.a.a.a.a.).

1.a.b.b.b) Let $\beta_{22}=0$. Then we are in case 1.a.a.b.a).

1.b) Let $\delta_{2}'\neq 0$. Then Example~7 with a scalar $\alpha = \delta_2'$ and multiplying elements $x_2,x_3$ by a scalar $\frac{1}{\delta_2'}$ allows us to assume that
\begin{gather*}
x_1 = e_{21} + (\alpha_{22}'+1)e_{11} +\alpha_{22}' ve_{12} + ve_{22},\\
x_2 = ve_{11} + \beta_{11}'e_{11} + (\beta_{22}'+\alpha_{22}'\beta_{22}'-1)e_{12} + \beta_{22}'e_{22} + \alpha_{22}'\beta_{11}'\beta_{22}'ve_{12} + \beta_{11}'\beta_{22}' ve_{22},\\
 x_3 = ve_{21} - e_{11} + (-\alpha_{22}'^2\beta_{22}'+2\alpha_{22}'-2\alpha_{22}'\beta_{22}'+1-\beta_{22}')e_{12} + (-\beta_{22}'+1-\alpha_{22}'\beta_{22}')e_{22} + \\ (-\alpha_{22}'\beta_{22}'+1-\beta_{22}') ve_{12} + (\beta_{11}'-\beta_{22}'-\beta_{11}'\beta_{22}') ve_{22},
\end{gather*}
where $\beta_{11}'=\frac{\beta_{11}}{\delta_2'}$, $\beta_{22}'=\frac{\beta_{22}}{\delta_2'}$. Example~5 with a scalar $\alpha = -1$, Example~1 with a scalar $\alpha = -(\beta_{11}'-\beta_{22}'-\beta_{11}'\beta_{22}')$ and Example~3 with a scalar $\alpha = -(2\alpha_{22}'+1)$ allow us to assume that
\begin{gather*}
x_1 = e_{21} + (\alpha_{22}'+1)e_{11},\\
x_2 = ve_{11} + \beta_{11}'e_{11} + \beta_{22}'(\alpha_{22}'+1)e_{12} + \beta_{22}'e_{22} + (\alpha_{22}'+1)\beta_{11}'\beta_{22}'ve_{12} + \beta_{11}'\beta_{22}' ve_{22},\\
 x_3 = ve_{21} -\beta_{22}'(\alpha_{22}'+1)^2e_{12} -\beta_{22}' (\alpha_{22}'+1)e_{22}.
\end{gather*}

1.b.a) Let $\alpha_{22}'+1\neq 0$. Then Example~13 with a scalar $\alpha = \frac{1}{\alpha_{22}'+1}$, multiplying elements $x_1,x_3$ by a scalar $\frac{1}{\alpha_{22}'+1}$ and Example~8 with a scalar $\beta_{11}'$ and considering element $x_2-\beta_{11}'x_1$ instead of element $x_2$ allow us to assume that
\begin{gather*}
x_1 = e_{21} + e_{11},\\
x_2 = ve_{11} + \beta_{22}'e_{12} + \beta_{22}'e_{22} + \beta_{11}'\beta_{22}'ve_{12},\\
x_3 = ve_{21} -\beta_{22}'e_{12} -\beta_{22}' e_{22} + \beta_{11}'\beta_{22}' ve_{22}.
\end{gather*}

1.b.a.a) Let $\beta_{22}'=0$. Then we get to case 1.a.a.a.a).

1.b.a.b) Let $\beta_{22}'\neq 0$. Then Example~7 with a scalar $\alpha = \beta_{22}'$, multiplication of elements $x_2,x_3$ by a scalar $\frac{1}{\beta_{22}'}$ and Example~1 with a scalar $\alpha = - \frac{\beta_{11}'}{\beta_{22'}}$ allows us to assume that we have fallen into case 1.a.a.a.b).

1.b.b) Let $\alpha_{22}'+1 = 0$.

1.b.b.a) Let $\beta_{11}'=0$. Then we are in case 1.a.b) with $\alpha_{22}'=0$

1.b.b.b) Let $\beta_{11}'\neq 0$, $\beta_{22}'=0$. Then the canonical involution allows us to assume that the coefficient of $e_{22}$ in $x_2$ is zero, which reduces us to case 1.b.b.a).

1.b.b.c) Let $\beta_{11}'\neq 0$, $\beta_{22}'\neq 0$. Then Example~7 with a scalar $\alpha = \beta_{11}'$ and the multiplication of elements $x_2,x_3$ by a scalar $\frac{1}{\beta_{11}'}$ allow us to assume that we have fallen into case 1.a.a.b.b).

2) Let $\alpha_{11}=0$, $\alpha_{22}\neq 0$. Then the canonical involution allows us to assume that the coefficient of $e_{11}$ in $x_1$ is nonzero, which reduces us to the first case.

3) Let $\alpha_{11}=\alpha_{22}=0$. Then
\begin{gather*}
 x_1 = e_{21} + \delta_2 ve_{22},\\
 x_2 = ve_{11} + \beta_{11}e_{11} -\delta_2e_{12} + \beta_{22}e_{22} + \beta_{11}\beta_{22} ve_{22},\\
 x_3 = ve_{21} -\delta_2 e_{11} + \delta_2e_{22} + \delta_2^2 ve_{12} + (\delta_2\beta_{11}-\delta_2\beta_{22}) ve_{22}.
\end{gather*}

3.a) Let $\delta_2 = 0$. Then we are in case 1.b.b).

3.b) Let $\delta_2\neq 0$. Then Example~7 with a scalar $\alpha = \delta_2$ and multiplication of elements $x_2,x_3$ by a scalar $\frac{1}{\delta_2}$ allow us to assume that
\begin{gather*}
x_1 = e_{21} + ve_{22},\\
x_2 = ve_{11} + \beta_{11}''e_{11} - e_{12} + \beta_{22}''e_{22} + \beta_{11}''\beta_{22}'' ve_{22},\\
x_3 = ve_{21} - e_{11} + e_{22} + ve_{12} + (\beta_{11}''-\beta_{22 }'') ve_{22},
\end{gather*}
where $\beta_{11}''=\frac{\beta_{11}}{\delta_2}$, $\beta_{22}''=\frac{\beta_ {22}}{\delta_2}$.

3.b.a) Let $\beta_{11}''=0$.

3.b.a.a) Let $\beta_{22}''=0$. Then Example~5 with a scalar $\alpha=-1$ allows us to assume that we are in case 1.a.a.b.a)

3.b.a.b) Let $\beta_{22}''\neq 0$. Then Example~12 with a scalar $\alpha = \beta_{22}''$, multiplying element $x_1$ by a scalar $\beta_{22}''$, multiplying element $x_2$ by a scalar $\frac{1}{ \beta_{22}''}$, Example~5 with a scalar $\alpha = -1$ and Example~1 with a scalar $\alpha = -1$ allow us to assume that we have fallen into case 1.a.b.b.a).

3.b.b) Let $\beta_{11}''\neq 0$, $\beta_{22}''=0$. Then the canonical involution allows us to assume that we are in case 3.b.a.

3.b.c) Let $\beta_{11}''\neq 0$, $\beta_{22}''\neq 0$. Then Example~12 with a scalar $\alpha = \beta_{22}''$, multiplication of element $x_1$ by a scalar $\beta_{22}''$ and multiplication of element $x_2$ by a scalar $\frac{1}{\beta_{22}''}$ allow us to assume that
\begin{gather*}
x_1 = e_{21} + ve_{22},\\
x_2 = ve_{11} + \delta e_{11} - e_{12} + e_{22} + \delta ve_{22},\\
x_3 = ve_{21} - e_{11} + e_{22} + ve_{12} + (\delta-1) ve_{22},
\end{gather*}
where $\delta = \frac{\beta_{11}''}{\beta_{22}''}$. Example~14 and the replacement $x_1\longleftrightarrow x_2$ allow us to assume that
\begin{gather*}
x_1 = e_{21} + e_{11} + \delta e_{12} + \delta e_{22} - ve_{22},\\
x_2 = ve_{11} + e_{12},\\
x_3 = ve_{21} + e_{11} + (\delta-1) e_{12} - e_{22} + ve_{12},
\end{gather*}
We are at the beginning of case 1) with scalars $\beta_{11}=\beta_{22}=0$, $\alpha_{22}'=\delta$, $\delta_2'=-1$. The lemma is proved.

\medskip In a four-dimensional case there are two different possibilities: a unital decomposition (one of subalgebras contains a unit of $\mathbb{O}$) and a non-unital decomposition. Let us start from unital one.

\medskip\textbf{Lemma 7.} {\it Let $\mathbb{O}=A_1\oplus A_2$ be a decomposition of the algebra of split octonions into a direct sum of two subalgebras, where $\dim A_1 = 4$ and $1\in A_1$. Let the field $F$ be a quadratically closed and have a characteristic different from 2. Then, up to automorphism and antiautomorphism, one of the following cases is possible:
\begin{enumerate}
\item $A_1 = Fe_{11}+Fe_{22}+Fve_{12}+Fve_{22}$, $A_2 = F(e_{11}+e_{12}) + F(e_{21}+e_{22}) + Fve_{11} + Fve_{21}$;
 \item $A_1 = F\cdot 1 + Fe_{12} + Fve_{12} + Fve_{22}$, $A_2 = Fe_{11} + Fe_{21} + Fve_{11} + Fve_{21}$.
\end{enumerate}
}

\textbf{Proof.} By the results of \cite{Octonions}, the non-unital four-dimensional subalgebra of the octonions is of the form $A_2 = \mathbb{O} x$ for some $x\in \mathbb{O}$, $x^2 = 0$.

Since $x^2 = 0$, then $\alpha_{11}+\alpha_{22}=0$. Moreover, since $n(x)=0$, then $n(A_2) = n(x\mathbb{O}) = 0$. Recall that for the element $y = c + vd$, where $c,d\in M_2(F)$, $n(y) = \det c - \det d$ holds.

Since the field is quadratically closed and has a characteristic different from 2, then, by virtue of the results of \cite{Octonions} up to automorphism and antiautomorphism, we can assume that we have one of the following cases.

1) $A_1 = M_2(F)$. Then in $A_2$ there is a basis of the following form:
\[
y_1 = ve_{11} + M_1,\quad
y_2 = ve_{12} + M_2,\quad
y_3 = ve_{21} + M_3,\quad
y_4 = ve_{22} + M_4,
\]
where $M_i\in M_2(F)$. Since $n(A_2)=0$, then $n(Fy_1+Fy_2) = 0$. This means that $\det(\xi M_1+\zeta M_2) = 0$ for any $\xi,\zeta\in F$. Since there is no two-dimensional subspace of matrices with zero determinant in the algebra of second-order matrices, then the matrices $M_1$ and $M_2$ are linearly dependent. Similarly, the matrices $M_1$ and $M_3$, $M_2$ and $M_4$, $M_3$ and $M_4$ are linearly dependent. Thus, we can assume that $M_i = \varepsilon_i M$ for some $\varepsilon_i\in F$ and $M\in M_2(F)$. But in this case $n(y_1+y_4) = \det(M_1+M_4) - 1 = -1$, it is a contradiction with the fact that $n(A_2) = 0$.

2) $A_1 = Fe_{11} + Fe_{22} + Fve_{12} + Fve_{22}$. Then $A_2$ has a basis of the following form:
\begin{gather*}
y_1 = e_{12} + \alpha e_{11} + \beta e_{22} + \gamma ve_{22} + \delta ve_{12} ,\\
y_2 = e_{21} + \alpha_1 e_{11} + \beta_1 e_{22} + \gamma_1 ve_{22} + \delta_1 ve_{12},\\
y_3 = ve_{21} + \alpha_2 e_{11} + \beta_2 e_{22} + \gamma_2 ve_{22} + \delta_2 ve_{12},\\
y_4 = ve_{11} + \alpha_3 e_{11} + \beta_3 e_{22} + \gamma_3 ve_{22} + \delta_3 ve_{12}.
\end{gather*}
Since $n(A_2)=0$, then $n( y_1 + y_2) = n(y_1)=n (y_2)=0$. Then we get $\alpha\beta = \alpha_1\beta_1 = 0$ and $(\alpha + \alpha_1)(\beta+\beta_1)-1 = 0$. If $\alpha = \alpha_1 = 0$ or $\beta = \beta_1 = 0$, then we have a contradiction. If $\alpha = \beta_1 = 0$, then we have $\alpha_1\beta= 1$. If $\alpha_1=\beta = 0$, then we have $\alpha\beta_1 = 1$. Up to the action of the canonical involution, we can assume that $\beta=\alpha_1 = 0$, $\alpha\beta_1 = 1 $.

We have $0 = n(y_3) = \alpha_2\beta_2 + \delta_2$, $0 = n(y_4) = \alpha_3\beta_3 - \gamma_3$.

Similarly, $0=n(y_1+y_4) = (\alpha+\alpha_3)\beta_3 - (\gamma+\alpha_3\beta_3) = \alpha\beta_3 - \gamma$, $n(y_2+ y_3) = \alpha_2(\frac {1}{\alpha} + \beta_2) +\delta_1 - \alpha_2\beta_2 = \frac{\alpha_2}{\alpha} + \delta_1$, $0 = n(y_2+y_4) = \alpha_3(\frac{1}{\alpha} + \beta_3) - (\gamma_1+\alpha_3\beta_3) = \frac{\alpha_3}{\alpha} - \gamma_1$ and $0 = n(y_1+y_3) = (\alpha+\alpha_2)\beta_2 + (\delta-\alpha_2\beta_2) = \alpha\beta_2+\delta$. Thus, $\gamma = \alpha\beta_3$, $\delta = -\alpha\beta_2$, $\alpha_3 = \gamma_1\alpha$ and $\alpha_2 = -\delta_1\alpha$.

Next,
$0 = n(y_3 + y_4) = (\alpha_2 + \alpha_3)(\beta_2 + \beta_3) - (\gamma_2+\gamma_3) + (\delta_2+\delta_3) = \alpha_3\beta_2 + \alpha_2\beta_3 - \gamma_2 +\delta_3 = \gamma_1\alpha\beta_2 - \delta_1\alpha\beta_3 - \gamma_2 + \delta_3$. Thus, $\delta_3 = \gamma_2 + \alpha(\delta_1\beta_3-\gamma_1\beta_2)$.

Since $y_1y_4 \in (\beta_3-\delta)e_{12} + A_1$, then $y_1y_4 = (\beta_3-\delta)y_1$. Comparison of the coefficients of $e_{11}$ in this equality yields the following equality:
\[\alpha(\alpha\beta_2-\alpha_1\gamma_1)=0,\]
from which we have $\beta_2 = \frac{\alpha_1\gamma_1}{\alpha}$.

Example~13 with a scalar $\alpha$, multiplication of $y_1$ by a scalar $\frac{1}{\alpha}$ and multiplication of $y_2,y_3$ by a scalar $\alpha$ allow us to consider that
\begin{gather*}
    y_1 = e_{12} + e_{11} + \beta_3 ve_{22} - \gamma_5 ve_{12},\\
    y_2 = e_{21} + e_{22} + \gamma_5 ve_{22} + \delta_4 ve_{12},\\
    y_3 = ve_{21} -\delta_4 e_{11} + \gamma_5 e_{22} + \gamma_4 ve_{22} +\delta_4\gamma_5 ve_{12},\\
    y_4 = ve_{11} + \gamma_5 e_{11} + \beta_3 e_{22} + \gamma_5\beta_3 ve_{22} + (\gamma_4 + \delta_4\beta_3-\gamma_5^2) ve_{12},
\end{gather*}
where $\delta_4 = \delta_1\alpha^2$, $\gamma_4 = \gamma_2\alpha$, $\gamma_5 = \gamma_1\alpha$. 

The automorphism from Example~3 with the scalar $-\delta_4$, a consideration of the element $y_3+\delta_4 y_1$ instead of the element $y_3$, the automorphism from Example~4 with the scalar $-(\delta_4\beta_3+\gamma_4)$, the automorphism from Example~2 with a scalar $\gamma_5$, a consideration of the element $y_3-\gamma_5 y_2$ instead of the element $y_3$, the automorphism from Example~1 with a scalar $\gamma_5^2$, the automorphism from Example~8 with a scalar $(\beta_3+\gamma_5)$ and a consideration of the element $y_4-(\beta_3+\gamma_5) y_2$ instead of the element $y_4$  allow us to consider that
\[ y_1 = e_{12} + e_{11},\quad y_2 = e_{21} + e_{22},\quad  y_3 = ve_{21},\quad
    y_4 = ve_{11}.\]

3) $A_1=F\cdot 1 + Fe_{12} + Fve_{12} + Fve_{22}$. Then in $A_2$ there is a basis of the following form:
\begin{gather*}
y_1 = e_{11} + \alpha 1 + \beta e_{12} + \gamma ve_{12} + \delta ve_{22},\\ y_2 = e_{21} + \alpha_1 1 + \beta_1 e_{12} + \gamma_1 ve_{12} + \delta_1 ve_{22},\\
y_3 = ve_{11} + \alpha_2 1 + \beta_2 e_{12 } + \gamma_2 ve_{12} + \delta_2 ve_{22},\\
y_4 = ve_{21} + \alpha_3 1 + \beta_3 e_{12} + \gamma_3 ve_{12} + \delta_3 ve_{22},
\end{gather*}

Since $n(A_2)=0$, then $0=n(y_1)=(\alpha+1)\alpha$, whence $\alpha \in\{0,-1\}$. Up to canonical involution and multiplying $y_1$ by $-1$, we can assume that $\alpha = 0$. Further, $0=n(y_2)=\ alpha_1^2 - \beta_1$, whence $\beta_1 = \alpha_1^2$. Further, $0=n(y_3)=\alpha_2^2 - \delta_2$, whence $\delta_2 = \alpha_2^2$. Next, $0=n(y_4) = \alpha_3^2 + \gamma_3$, whence $\gamma_3 = -\alpha_3^2$. Next, $0=n(y_1+y_2)=(\alpha_1+1)(\alpha_1 )-\beta$, whence $\beta = (\alpha_1+1)\alpha_1$. Further, $0=n(y_1+y_3) = (\alpha_2+1)\alpha_2-\delta-\delta_2 = \alpha_2- \delta$, whence $\delta = \alpha_2$. Further, $0=n(y_1+y_4) = (\alpha_3+1)\alpha_3+\gamma+\gamma_3 = \alpha_3 + \gamma$, whence $\gamma = - \alpha_3$. Next, $0=n(y_2+y_3) = (\alpha_1+\alpha_2)^2-(\beta_1+\beta_2)-\delta_1-\delta_2 = 2\alpha_1\alpha_2-\beta_2-\delta_1$. Next, $0=n(y_2+y_4) = (\alpha_1+\alpha_3)^2-\beta_1-\beta_3+\gamma_1+\gamma_3 = 2\alpha_1\alpha_3-\beta_3-\gamma_1$. Next, $0=n(y_3+y_4) = (\alpha_2+\alpha_3)^2-\delta_2-\delta_3+\gamma_2+\gamma_3 = 2\alpha_2\alpha_3-\delta_3+\gamma_2$.

Since $y_2y_1\in -\alpha_1 ^2 e_{11} + e_{21} +A_1$, then $y_2y_1=-\alpha_1^2 y_1 + y_2$. Comparing the coefficients of the unit in this equality, we get $\alpha_1^2=0$, whence $ \alpha_1 = 0$. Then, comparing the coefficients of $ve_{22}$, we get $\alpha_3=\delta_1$. Further, $y_2y_3\in \delta_1 e_{11}+\alpha_2 e_{21}+ve_{21}+A_1$, so that $y_2y_3=\delta_1 y_1 + \alpha_2 y_2 +y_4$. But direct calculations show that $y_2y_3-\delta_1y_1 - \alpha_2 y_2 -y_4 = -2\gamma_1 e_{12}$, whence $\gamma_1 = 0$. Example~1 with scalar $-\gamma_2$ allows us to assume that $\gamma_2 = 0$. Example~5 with a scalar $-\delta_1$ and consideration of the element $y_4-2\delta_1 y_1$ instead of the element $y_4$ allow us to assume that $\delta_1 = 0$. Example~16 with a scalar $-\alpha_2$ and consideration of the element $y_3-2\alpha_2 y_1$ instead of the element $y_3$ allow us to assume that $\alpha_2 = 0$. Thus, we get
\[
y_1 = e_{11} ,\quad
y_2 = e_{21},\quad
y_3 = ve_{11},\quad
y_4 = ve_{21}.
\]
The lemma is proved.

\medskip Let us note that a condition of quadratically closure is important.

\medskip\textbf{Example 19.} \textit{Let $\mathbb{O}$ be a split Cayley-Dickson algebra over a field $\mathbb{R}$. We have a decomposition $\mathbb{O}=A_1\oplus A_2$, where $A_1 = F\cdot 1 + F(e_{12}-e_{21}) + Fve_{12} + Fve_{22}$, $A_2 = Fe_{11} + Fe_{21} + Fve_{11} + Fve_{21}$. An subalgebra $A_1$ is not isomorphic to the unital subalgebras $C_1 = Fe_{11}+Fe_{22}+Fve_{11}+Fve_{12}$ and $C_2 = F\cdot 1 + Fe_{12} + Fve_{12} + Fve_{22}$ in Lemma~7. First, $A_1$ is not isomorphic to $C_1$ because its radicals have different dimensions. For $C_2$ we have $(e_{12}-e_{21})^2 = -1$. If $A_1\simeq C_2$ then $C_2$ contains an element $x$, such that $x^2 = -1$. Let $x = \alpha\cdot 1 + \beta e_{12} + \gamma ve_{12} + \delta ve_{22}$. Then $x^2 = \alpha^2 \cdot 1 + 2\alpha\beta e_{12} + 2\alpha\gamma ve_{12} + 2\alpha\delta ve_{22}$. It means that $x^2 = -1$ iff $\alpha^2 = -1$ and $\beta = \gamma = \delta = 0$. So, $A_1$ is not isomorphic to $C_2$. 
}

\medskip Let us see a non-unital case of decomposition on four-dimensional subalgebras.

\medskip\textbf{Lemma 8.} {\it Let $\mathbb{O} = A_1 \oplus A_2$ be the decomposition of the algebra of split octonions into a direct sum of two non-unital subalgebras, where $\dim A_1 = 4$. Then, up to automorphism and antiautomorphism, we can assume that $A_1 = Fe_{11} + Fe_{12} + Fve_{11} + Fve_{12}$ and $A_2 = Fe_{21} + Fe_{22} + Fve_{21} + Fve_{22}.$
}

\textbf{Proof.} Let $\mathbb{O} = A_1 \oplus A_2$, where $\dim A_1 = \dim A_2 = 4$, $A_1$ and $A_2$ are subalgebras of $\mathbb{O}$. Then, up to automorphism and antiautomorphism, we can assume that $A_1= Fe_{11} + Fe_{12} + Fve_{11} + Fve_{12}$ (see \cite{Octonions}). Then $A_2$ has a basis of the following form:
\begin{gather*}
y_1 = e_{21} + \alpha e_{11} + \beta e_{12} + \gamma ve_{11} + \delta ve_{12},\\
y_2 = e_{22} + \alpha_1 e_{11} + \beta_1 e_{12} + \gamma_1 ve_{11} + \delta_1 ve_{12},\\
y_3 = ve_{21} + \alpha_2 e_{11} + \beta_2 e_{12} + \gamma_2 ve_{11} + \delta_2 ve_{12},\\
y_4 = ve_{22} + \alpha_3 e_{11} + \beta_3 e_{12} + \gamma_3 ve_{11} + \delta_3 ve_{12}.
\end{gather*}

Note that $y_2^2 \in e_{22} + A_1$. This means that $y_2^2 = y_2$. From this equality, comparing the coefficients of $e_{11}$, $e_{12}$, $ve_{11}$, $ve_{12}$, we get $\alpha_{1}^2 = \alpha_1$, $\alpha_1\beta_1 = \alpha_1\gamma_1=\alpha_1\delta_1 = 0$.

Let us consider some cases.

1) $\alpha_{1} = 1$, then from the above $\beta_1=\gamma_1=\delta_1 = 0$. In particular, $y_2 = 1$, it is a contradiction.

2) $\alpha_1 = 0$. Then from $y_1^2 \in \alpha e_{21} + \beta e_{22} + A_1$ we get $y_1^2 = \alpha y_1 + \beta y_2$. Comparing the coefficients of $e_{11}$, we get $\beta = 0$. Similarly, from $y_3^2$ and $y_4^2$ we get $\delta_2 = \gamma_3 = 0$.

Next, $y_3y_4 \in -e_{21}+\gamma_2 e_{22} + \alpha_2 ve_{12} + A_1$, whence $y_3y_4 = -y_1 +\gamma_2 y_2 + \alpha_2 y_4$. Comparing the coefficients of $e_{11}$ in this equality, we get $\alpha = \delta_3$. Similarly, from $y_4y_3$ we get $\alpha = \gamma_2$. Similarly, from $y_4y_2$ we get $\gamma_1 = \alpha_3$, whence, comparing the coefficients of $ve_{12}$ we get $\alpha = \beta_1$. Similarly, from $y_3y_1$ we get $\delta = \beta_2$. Similarly, from $y_3y_2$ we get $\delta_1 = -\alpha_2$. Then from $y_3y_4$ we get $\gamma = -\beta_3$. Example~3 with a scalar $-\delta$ allows us to assume that $\delta = 0$. Example~4 with a scalar $-\gamma$ allows us to assume that $\gamma = 0$.

2.a) Let $\alpha = 0$.

2.a.a) Let $\delta_1=\gamma_1 = 0$. Then
\[ y_1 = e_{21},\quad y_2 = e_{22},\quad y_3 = ve_{21},\quad y_4 = ve_{22}.\]

2.a.b) Let $\delta_1 \neq 0$, $\gamma_1 = 0$. Then the automorphism from Example~13 with a scalar $\frac{1}{\delta_1}$, multiplication of elements $y_1,y_3$ by $\frac{1}{\delta_1}$,
the automorphism from Example~5 with a scalar $-1$, a consideration of an element $y_1+y_4$ instead of an element $y_1$ and a consideration of an element $y_3+y_2$ instead of an element $y_3$ allow us to assume that
\[ y_1 = e_{21},\quad y_2 = e_{22},\quad y_3 = ve_{21},\quad y_4 = ve_{22}.\]

2.a.c) Let $\gamma_1\neq 0$, $\delta_1 = 0$. Then the automorphism from Example~9 and the sign change of the scalar $\gamma_1$ leads us to case 2.a.b.

2.a.d) Let $\gamma_1\neq 0$, $\delta_1 \neq 0$. Then the automorphism from Example~13 with the scalar $\frac{1}{\delta_1}$, the automorphism from Example~12 with the scalar $\frac{1}{\gamma_1}$, the multiplication of the element $y_1$ by $\frac{1}{\gamma_1\delta_1}$, the multiplication of the element $y_3$ by $\frac{1}{\delta_1}$, the multiplication of the element $y_4$ by $\frac{1}{\gamma_1}$ allow us to assume that
\[
    y_1 = e_{21},\quad
    y_2 = e_{22} + ve_{11} + ve_{12},\quad
    y_3 = ve_{21} - e_{11},\quad
    y_4 = ve_{22} + e_{11}.
\]
Then the Example~18 allows us to assume that 
\[ y_1 = e_{21},\quad y_2 = e_{22},\quad y_3 = ve_{21},\quad y_4 = ve_{22}.\]

2.b) Let $\alpha \neq 0$. Then the automorphism from Example~13 with scalar $\frac{1}{\alpha}$ and the multiplication of an element $y_3$ by a scalar $\frac{1}{\alpha}$ allow us to assume that
\begin{gather*}
y_1 = e_{21} + e_{11},\\
y_2 = e_{22} + e_{12} + \gamma_1 ve_{11} + \delta_1' ve_{12},\\
y_3 = ve_{21} - \delta_1' e_{11} + ve_{11},\\
y_4 = ve_{22} + \gamma_1 e_{11} + ve_{12},
\end{gather*}
where $\delta_1'=\frac{\delta_1}{\alpha}$. Example~6 with scalar $1$ and a consideration of an element $y_1-y_2$ instead of an element $y_1$ allows us to assume that
\begin{gather*}
y_1 = e_{21} - \gamma_1 ve_{11} - \delta_1' ve_{12},\\
y_2 = e_{22} + \gamma_1 ve_{11} + \delta_1' ve_{12},\\
y_3 = ve_{21} - \delta_1' e_{11} - \delta_1'e_{12},\\
y_4 = ve_{22} + \gamma_1 e_{11} + \gamma_1 e_{12},
\end{gather*}
Then Example~3 with a scalar $\delta_1'$ and Example~4 with a scalar $\gamma_1$ allow us to assume that we have a case 2.a. The lemma is proved.

\medskip Lemmas 5--8 imply the following statement. 

\medskip\textbf{Theorem 2.} {\it Let $\mathbb{O} = A_1 \oplus A_2$ be a decomposition of the algebra of split octonions into a direct sum of two subalgebras over a quadratically closed field with characteristic $\neq 2$. Then, up to automorphism and antiautomorphism, one of the following cases holds:
\begin{enumerate}
\item $A_1 = M_2(F)+Fve_{12}+Fve_{22}$, $A_2 = Fve_{11} + Fve_{21}$;
\item $A_1 = M_2(F)+Fve_{12}+Fve_{22}$, $A_2 = Fve_{11} + F(ve_{21}+e_{22})$;
\item $A_1 = Fe_{11} + Fe_{12} + Fe_{22} + Fve_{12} + Fve_{22}$, $A_2 = Fe_{21} + Fve_{11} + Fve_{21}$;
\item $A_1 = Fe_{11} + Fe_{12} + Fe_{22} + Fve_{12} + Fve_{22}$, $A_2 = F(e_{21}+e_{11}) + Fve_{11} + Fve_{21}$;
\item $A_1 = Fe_{11}+Fe_{22}+Fve_{12}+Fve_{22}$, $A_2 = F(e_{11}+e_{12}) + F(e_{21}+e_{22}) + Fve_{11} + Fve_{21}$;
\item $A_1 = F\cdot 1 + Fe_{12} + Fve_{12} + Fve_{22}$, $A_2 = Fe_{11} + Fe_{21} + Fve_{11} + Fve_{21}$;
\item $A_1 = Fe_{11} + Fe_{12} + Fve_{11} + Fve_{12}$, $A_2 = Fe_{21} + Fe_{22} + Fve_{21} +  Fve_{22}$.
\end{enumerate}}

\medskip This theorem completes the description of Rota-Baxter operators on split composition algebras (see \cite{BGP}, \cite{RBZ},  \cite{Unital}). It is easy to see that all 7~decompositions in Theorem~2 are non-isomorphic to each other.

   \section{Acknowlegements}

The author expresses gratitude to V.~Yu.~Gubarev for useful discussions.

The study was supported by a~grant from the Russian Science Foundation №~23-71-10005, https://rscf.ru/project/23-71-10005/

\end{document}